\documentclass[pdflatex,sn-mathphys-num]{sn-jnl}

\usepackage{graphicx}%
\usepackage{multirow}%
\usepackage{amsmath,amssymb,amsfonts}%

\usepackage{amsthm}%
\usepackage{mathrsfs}%
\usepackage[title]{appendix}%
\usepackage{xcolor}%
\usepackage{textcomp}%
\usepackage{manyfoot}%
\usepackage{booktabs}%
\usepackage{algorithm}%
\usepackage{algorithmicx}%
\usepackage{algpseudocode}%
\usepackage{listings}%
\usepackage{siunitx}
\usepackage{longtable,tabularx}
\usepackage{comment}
\usepackage{upgreek}
\usepackage{kotex}
\usepackage{makecell}
\usepackage{setspace}
\usepackage{hyperref}
\usepackage{cleveref}
\usepackage{mathtools}
\usepackage{float}
\usepackage{enumitem}
\usepackage{subcaption}

\newtheorem{theorem}{Theorem}%
\newtheorem{proposition}[theorem]{Proposition}%
\newtheorem{lemma}{Lemma}
\newtheorem{assumption}{Assumption}

\newtheorem{remark}{Remark}%

\newtheorem{definition}{Definition}%

\begin{document}

\title[Feedback Integrators Revisited]{Feedback Integrators: Non-Asymptotic Invariance for One-Step Methods and Gain Selection under Euler Discretization}


\author[1]{\fnm{Juho} \sur{Bae}}\email{juhobae@mit.edu}

\author*[1]{\fnm{Dong Eui} \sur{Chang}}\email{dechang@kaist.ac.kr}


\affil[1]{\orgdiv{School of Electrical Engineering}, \orgname{Korea Advanced Institute of Science and Technology}, \orgaddress{\city{Daejeon}, \postcode{34141}, \country{Republic of Korea}}}

\abstract{For dynamical systems evolving on a manifold and admitting first integrals, standard one-step numerical methods generally cause the discrete trajectory to drift off the manifold and the numerical values of the first integrals to deviate from their prescribed values. Feedback integrators address this by extending the dynamics to an ambient Euclidean space and adding a feedback term that drives the numerical trajectory toward the set satisfying both the manifold constraint and the prescribed values of the first integrals. Existing theory, however, has two limitations: it remains asymptotic, guaranteeing only eventual entrance into an attractor containing the desired set, and it does not explain how the feedback gain should be chosen. In this paper, we first close the former gap for general one-step methods by proving positive invariance of arbitrarily small sublevel neighborhoods of the feedback Lyapunov function for sufficiently small step sizes. We then specialize to Euler discretization and analyze how the feedback gain enters the Taylor-based error bound. In this setting, we characterize a range of scaled gains that guarantee positive invariance for sufficiently small step sizes and identify the scaling that minimizes the Taylor-based upper bound. We further propose adaptive gain-selection rules under Euler discretization, including both stepwise and periodically updated variants, and establish corresponding boundedness guarantees for the resulting discrete trajectories. These results identify Euler discretization as the first setting in which gain selection for feedback integrators closes in explicit form, whereas extensions to general higher-order one-step methods remain genuinely method-dependent. Numerical experiments on free rigid body motion in $\operatorname{SO}(3)$, the Kepler problem, and a perturbed Kepler problem with rotational symmetry support the analysis.
}


\maketitle

\section{Introduction}\label{sec:intro}
Consider a dynamical system 
\begin{equation} \label{eq:original_system}
	\dot{x} = f(x), \quad x(0) = x_I,
\end{equation}
defined on a smooth manifold $M$, where the dynamics $f:M \rightarrow TM$ is $C^1$ and solutions are assumed to exist globally. Suppose that the system admits first integrals $f_j: M \rightarrow \mathbb{R}$, $j \in \{ 1, 2, \dots \ell \}$. 

Throughout the feedback-integrator construction, we assume that $M$ is embedded in $\mathbb R^n$, and that $f$ and the first integrals $f_j$ admit smooth extensions to an open neighborhood $U\subset\mathbb R^n$ of $M$. In this ambient representation, an ordinary one-step method with step size $h$ may be written as
\begin{equation} \label{eq:original_discrete_system}
	x_{k+1} = x_k + h f_h(x_k), \quad x_0 = x_I.
\end{equation}
There is, however, no guarantee that the discrete trajectory remains on $M$ or that the numerical values of first integrals are preserved. The truncation error introduced during discretization may cause the discrete state to leave $M$ and the numerical values of the first integrals to vary. There have been significant efforts to mitigate, if not resolve, such inconsistency introduced during numerical integration~\cite{Hairer2006}. 

Many structure-preserving methods are tailored to particular geometric structures or require modifying the numerical method itself, for instance through projection, splitting, symplectic, or Lie-group constructions~\cite{Hairer2006}. The \emph{Feedback Integrator} framework~\cite{chang2016feedback} takes a different approach: it modifies the dynamics rather than the integrator, allowing ordinary one-step methods to be applied to the resulting surrogate system. The key idea is as follows. The challenge of preserving the manifold structure and the first integrals is viewed as a stabilization problem in the ambient Euclidean space $\mathbb{R}^n$ in which $M$ is embedded. The target set $\Lambda \subset \mathbb{R}^n$ consists of the points satisfying both the manifold constraint and the prescribed values of the first integrals. The feedback term that is designed to push the state toward $\Lambda$ outside $\Lambda$ is then added to the extended continuous dynamics, so that $\Lambda$ becomes an attracting set for the closed-loop system in continuous time. The original paper~\cite{chang2016feedback} then proves practical asymptotic stability of the discretized closed-loop system with respect to $\Lambda$ by means of attractor theory of ODEs~\cite{Kloeden1986}.

To elaborate the concepts, we define the target set $\Lambda$, which is the set satisfying both the manifold constraint and the prescribed values of the first integrals.
\begin{equation}
	\Lambda \coloneqq \left\{ x \in U : x \in M, \enspace f_j(x) = f_j(x_I), \enspace j = 1, 2, \dots \ell \right\}
\end{equation}
Now assume there exists an analytic function $V: U \rightarrow \mathbb{R}_{\geq 0}$ with $V^{-1}(0) = \Lambda$ such that 
\begin{assumption}[Eligibility conditions for $V$] \label{assumption:01}
	\leavevmode
	\begin{enumerate}[label=(A\arabic*), ref=A\arabic*, leftmargin=*]
		\item \label{assumption:01:A1}
		$\langle \nabla V(x), f(x) \rangle = 0$ for all $x \in U$,

		\item \label{assumption:01:A2}
		there exists $\nu > 0$ such that $V^{-1}([0, \nu]) \subseteq U$ is compact,

		\item \label{assumption:01:A3}
		all critical points of $V$ in $V^{-1}([0, \nu])$ are in $V^{-1}(0)$.
	\end{enumerate}
\end{assumption}
Let us further denote the sublevel set of $V$ at $\varepsilon > 0$ as $U^\varepsilon \coloneqq V^{-1}([0, \varepsilon])$. We note that these assumptions and the {\L}ojasiewicz inequality further imply the existence of class--$\mathcal{K}$ functions $M, m: [0,\nu] \rightarrow \mathbb{R}_{\geq 0}$ such that for all $v \in [0, \nu]$, 
\begin{equation} \label{eq:gradient_bounds}
m(v) \leq \inf_{\substack{x\in U^\nu\\ V(x)\geq v}} |\nabla V(x)|, \quad \sup\limits_{\substack{x\in U^\nu\\ V(x)\leq v}} |\nabla V(x)| \leq M(v).
\end{equation}

Under Assumption~\ref{assumption:01}, the following \emph{surrogate} system is considered on $U$ with a \emph{feedback} term $- \alpha \nabla V(x)$ and gain $\alpha > 0$:
\begin{equation} \label{eq:conti_system}
	\dot{x} = Y(x) \coloneqq f(x) - \alpha \nabla V(x).
\end{equation}
System~\eqref{eq:conti_system} is a surrogate of system~\eqref{eq:original_system} in a sense that the two vector fields coincide on $\Lambda$. (Recall that $\nabla V = 0$ on $\Lambda$.) Consequently, any solution of~\eqref{eq:conti_system} with an initial value on $\Lambda$ coincides with that of~\eqref{eq:original_system} and thus enjoys global existence property. Note that one cannot yet state anything about global existence of solutions to~\eqref{eq:original_system} and~\eqref{eq:conti_system} on $U$. For any solution of~\eqref{eq:conti_system} with initial value in $U^\nu$, the identity
\begin{equation} \label{eq:lyapunov_decrease}
	\frac{d}{dt}V(x(t)) = \langle \nabla V(x(t)), f(x(t)) - \alpha \nabla V(x(t)) \rangle = -\alpha \left| \nabla V(x(t)) \right|^2 \leq 0,
\end{equation}
holds as long as the solution exists, where $|\cdot|$ denotes the Euclidean norm of vectors. Hence, $U^\nu$ is positively invariant. Since this set is compact and contained in $U$, solutions starting in $U^\nu$ cannot blow up in finite time and exist globally. Moreover, by Assumption~\ref{assumption:01}, the only critical points of $V$ in $U^\nu$ lie in $\Lambda$ and therefore, the surrogate system is asymptotically stable with respect to $\Lambda$ on $U^\nu$.

The key idea of feedback integrators is to numerically integrate the surrogate system~\eqref{eq:conti_system}, whose continuous-time dynamics makes $\Lambda$ attracting. The corresponding discrete-time stability, however, is not automatic. Technically, the discretized system considered in numerical integration may fail to have the same attractor set $\Lambda$, or it may even be unstable. To address this, the original paper of Chang~\cite{chang2016feedback} presents asymptotic performance guarantee with the following theorem.
\begin{theorem}[Theorem 5.2,~\cite{chang2016feedback}] \label{thm:chang2016}
	Suppose that a one--step method of order $p$ is applied to numerically integrate~\eqref{eq:conti_system}, denoted as follows.
\begin{equation} \label{eq:feedback_discrete_system}
	x_{k+1} = x_k + h Y_h(x_k), \quad x_0 = x_I.
\end{equation}
Suppose that the vector field $f$ is $C^{p}$ and Assumption~\ref{assumption:01} holds. Then there exists $h_{0}>0$ such that for all $h \in \left(0, h_0\right)$, the discrete system~\eqref{eq:feedback_discrete_system} has a compact, uniformly asymptotically stable set $\Lambda_{h}$ which contains $\Lambda$, and $\Lambda_{h} \rightarrow \Lambda$ as $h \rightarrow 0^{+}$ with respect to Hausdorff metric. Moreover, there exist (i) a bounded open set $U_{0}$ independent of $h$ and contains $\Lambda_{h}$, and (ii) a time $T_{0}(h)=A + B p\log \frac{1}{h}$ where $A$ and $B$ are constants, such that the sequence generated by~\eqref{eq:feedback_discrete_system} with $x_I \in U_0$ and $h \in (0, h_0)$ satisfies $x_{k} \in \Lambda_{h}$ for all $kh \geq T_{0}(h)$.
\end{theorem}
The notion of feedback integrator was further extended to nonholonomic systems~\cite{chang2019feedback} and systems with holonomic constraints~\cite{chang2022feedback}. A central advantage of the feedback integrator framework is its scheme-level generality: once the surrogate dynamics is constructed, it can be integrated by any one-step method without modifying the numerical scheme itself. This generality is also reflected in Theorem~\ref{thm:chang2016}, which is stated for arbitrary one-step methods of a given order. At the same time, as we discuss below, the same generality makes gain selection method-dependent once one asks for non-asymptotic analysis.

Theorem~\ref{thm:chang2016} leaves two coupled limitations in the existing theory. First, its guarantee remains asymptotic: it is formulated through an entry time $T_0(h)$, which diverges as $h\to0^+$. Hence, even though the initial condition satisfies $x_I\in\Lambda$, the theorem does not control the discrete trajectory over the whole integration horizon, but only after an asymptotic transient. This is particularly restrictive for gain selection, since the effect of the feedback gain $\alpha$ must be assessed from the first step onward, not only after the trajectory has entered the $h$-dependent attractor $\Lambda_h$.

Second, the original framework does not provide a criterion for choosing the gain in connection with such entire-horizon bounds. In particular, it does not identify which gain scalings are compatible with positive invariance in the small-step regime, nor does it provide a basis for deciding when one gain is preferable to another. This is not merely a normalization ambiguity, although scaling $V$ by a positive constant preserves Assumption~\ref{assumption:01}. Rather, the gain enters the discrete dynamics through the one-step discretization of the surrogate vector field $f-\alpha\nabla V$, and this dependence is determined by the chosen numerical method. From a control-theoretic viewpoint, this is natural: once the feedback term is implemented through a numerical method, stability and performance are governed by the resulting discrete-time closed-loop map, not by the continuous feedback term alone. At the level of general one-step methods, this method dependence appears through the local error expansion of the surrogate vector field. Consequently, gain selection should generally be formulated at the level of a specified discretization, rather than as a method-agnostic rule for arbitrary one-step methods. This motivates treating the problem first in Euler discretization, the canonical explicit case in which the additional method-dependent remainder terms are absent, and the gain dependence is exposed in closed form.

In this paper, we address these limitations in two stages. We first establish positive invariance of $U^\varepsilon$ under feedback integration for general one-step methods, thereby providing a non-asymptotic-in-time bound on the discrete trajectory. We then specialize to Euler discretization and obtain a complete small-step gain-selection theory in this canonical explicit setting: in terms of the scaled gain $\beta=h\alpha$, we characterize a range of $\beta$ that guarantees positive invariance for all sufficiently small $h$, identify the scaling that minimizes the Taylor-based upper bound, and develop stepwise and periodically updated adaptive gain-selection rules with corresponding boundedness guarantees.
\begin{remark}
	In~\cite{chang2016feedback}, candidate Lyapunov functions are constructed in sums-of-squares form and the eligibility conditions in Assumption~\ref{assumption:01} are verified under suitable assumptions. We use the same construction in the numerical demonstrations. The non-asymptotic invariance result below relies only on Assumption~\ref{assumption:01}; additional regularity assumptions are stated explicitly when needed for the gain-selection results.
\end{remark}
\section{Main Results}\label{sec:main}
\subsection{Preservation of First Integrals} \label{subsec:positive_invariance}
We first address the non-asymptotic gap in Theorem~\ref{thm:chang2016}. For every $\varepsilon \in (0,\nu)$, the next theorem establishes positive invariance of the sublevel set $U^\varepsilon = V^{-1}([0,\varepsilon])$ under the discrete system~\eqref{eq:feedback_discrete_system} for sufficiently small $h$. This fixed gain result will serve as the starting point for the gain-selection analysis in Sections~\ref{subsec:gain_selection} and~\ref{subsec:adaptive_gain}.
\begin{theorem} \label{thm:pos_invariance}
	For any $\varepsilon \in (0, \nu)$, there exists $h_0 > 0$ such that for all $0 < h < h_0$, $U^\varepsilon$ is positively invariant under the discrete system~\eqref{eq:feedback_discrete_system}.
\end{theorem}
\begin{proof}
	Let $x(\cdot;\xi)$ denote the solution of~\eqref{eq:conti_system} with initial value $\xi$. It suffices to show that $V(x_k)\leq \varepsilon \Longrightarrow V(x_{k+1})\leq \varepsilon$ for all sufficiently small $h$. Fix $c\in (\varepsilon,\nu)$. Since $U^\varepsilon$ is a compact subset of the open set $U$, there exists $r>0$ such that $(U^\varepsilon)_r \coloneqq \{x\in \mathbb{R}^n : \operatorname{dist}(x,U^\varepsilon)\leq r\} \subseteq U$. 

	From the Lyapunov decrease condition~\eqref{eq:lyapunov_decrease}, $V$ is nonincreasing along the solution of~\eqref{eq:conti_system}. Thus for $x_k\in U^\varepsilon$,
	\begin{equation}
		V(x(h;x_k)) \leq V(x_k)\leq \varepsilon,
	\end{equation}
	and $x(h;x_k)\in U^\varepsilon$. Since the one--step method has order $p$,
	there exist $h_{\mathrm{loc}},C_{\mathrm{loc}}>0$ such that
	\begin{equation}
		e \coloneqq x_{k+1}-x(h;x_k), \qquad |e|\leq C_{\mathrm{loc}} h^{p+1},
	\end{equation}
	for all $x_k\in U^\varepsilon$ and $h\in (0,h_{\mathrm{loc}}]$. Shrinking
	$h_{\mathrm{loc}}$ if necessary, we may assume $|e|<r$. Then $x_{k+1}\in B(x(h;x_k),r)\subseteq U$, and since the ball is convex, the segment joining $x(h;x_k)$ and $x_{k+1}$
	is contained in $U$. Because $V\in C^2(U)$, there exists $L_V>0$ such that $\left\|\nabla^2 V(x)\right\|\leq L_V$ for all $x\in (U^\varepsilon)_r$. Then Taylor's theorem gives
	\begin{equation}
		V(x_{k+1}) \leq V(x(h;x_k)) + \nabla V(x(h;x_k))\cdot e + \frac{L_V}{2}|e|^2.
	\end{equation}
	Since $\sup\limits_{x\in (U^\varepsilon)_r} |\nabla V(x)|<\infty$, it follows that
	\begin{equation} \label{eq:V_taylor1}
		V(x_{k+1}) \leq V(x(h;x_k)) + C h^{p+1}
	\end{equation}
	for some constant $C>0$.

	Set $\delta \coloneqq V(x_k),$ and $\delta' \coloneqq V(x(h;x_k))$. Since $V$ is nonincreasing along $x(t;x_k)$, we have $\delta' \leq V(x(t;x_k)) \leq \delta$ for $t\in [0,h]$. Therefore, we can introduce the following bound using the class--$\mathcal{K}$ functions $m$ and $M$ in~\eqref{eq:gradient_bounds} as
	\begin{equation}
		m(\delta') \leq |\nabla V(x(t;x_k))| \leq M(\delta),
		\qquad t\in [0,h].
	\end{equation}
	Using $\dot V = -\alpha |\nabla V|^2$, $\dot{V}(x)$ can be bounded as 
	\begin{equation}
		-\alpha M(\delta)^2 \leq \dot{V}(x) \leq -\alpha m\left( \delta' \right)^2,
	\end{equation}
	Integrating over $[0,h]$, we obtain
	\begin{equation}
		\delta - h\alpha M(\delta)^2 \leq \delta' \leq \delta - h\alpha m(\delta')^2.
		\label{eq:delta_prime_bound}
	\end{equation}

	Now, consider two cases $\delta \leq \varepsilon/2$ and $\delta \in (\varepsilon/2,\varepsilon]$ separately. If $\delta \leq \varepsilon/2$, then \eqref{eq:V_taylor1} gives
	\begin{equation}
		V(x_{k+1}) \leq \delta' + Ch^{p+1} \leq \frac{\varepsilon}{2} + Ch^{p+1}.
	\end{equation}
	Hence, $V(x_{k+1})\leq \varepsilon$ whenever $h \leq h_1 \coloneqq \min\left\{ h_{\mathrm{loc}}, \left(\frac{\varepsilon}{2C}\right)^{1/(p+1)} \right\}$.

	Assume next that $\delta \in (\varepsilon/2,\varepsilon]$. Let $h_2 \coloneqq \frac{\varepsilon}{4\alpha M(\varepsilon)^2}$. For $h<h_2$, \eqref{eq:delta_prime_bound} yields
	\begin{equation}
		\delta' \geq \delta - h\alpha M(\delta)^2 \geq \frac{\varepsilon}{2} - h\alpha M(\varepsilon)^2 > \frac{\varepsilon}{4}.
	\end{equation}
	Since $m$ is increasing,
	\begin{equation}
		m(\delta') \geq m\!\left(\frac{\varepsilon}{4}\right).
	\end{equation}
	Substituting this into the upper bound in \eqref{eq:delta_prime_bound}, we obtain
	\begin{equation}
		\delta'
		\leq \delta - h\alpha m(\delta')^2
		\leq \varepsilon - h\alpha m\!\left(\frac{\varepsilon}{4}\right)^2.
	\end{equation}
	Combining this with \eqref{eq:V_taylor1},
	\begin{equation}
		V(x_{k+1})
		\leq
		\varepsilon - h\alpha m\!\left(\frac{\varepsilon}{4}\right)^2 + Ch^{p+1}.
	\end{equation}
	Thus $V(x_{k+1})\leq \varepsilon$ whenever $h \leq h_3 \coloneqq \min\!\left\{ h_{\mathrm{loc}}, h_2, \left(\frac{\alpha\, m(\varepsilon/4)^2}{C}\right)^{1/p} \right\}$.

	Finally, set $h_0 \coloneqq \min\{h_1, h_3\}$. Then $V(x_k)\leq \varepsilon$ implies $V(x_{k+1})\leq \varepsilon$ for all $h\in (0,h_0)$. Therefore, $U^\varepsilon$ is positively invariant under~\eqref{eq:feedback_discrete_system}.
\end{proof}
\subsection{Fixed Gain Selection under Euler Discretization} \label{subsec:gain_selection}
We now focus on Euler discretization and examine how the gain enters the one--step change of $V$. Theorem~\ref{thm:pos_invariance} guarantees positive invariance of $U^\varepsilon$ for sufficiently small $h$ at fixed gain. To compare different gains, it is therefore enough to analyze a single Euler step issued from $U^\varepsilon$ and the Hessian of $V$ along the corresponding one--step segment.
\begin{definition} \label{def:01}
Fix $\varepsilon \in (0,\nu)$ and let $L \coloneqq \sup\limits_{x\in U^\varepsilon}\left\|\nabla^2 V(x)\right\|$ where $\|\cdot\|$ denotes the matrix 2-norm. For $h>0$ and $\alpha>0$, define
\begin{equation}
	U^\varepsilon_{h,\alpha} \coloneqq \left\{ x+t\,h\bigl(f(x)-\alpha \nabla V(x)\bigr) : x\in U^\varepsilon,\ t\in[0,1] \right\}.
\end{equation}
Whenever $U^\varepsilon_{h,\alpha}\subseteq U$, define
	\begin{equation}
L(h,\alpha) \coloneqq \sup\limits_{x\in U^\varepsilon_{h,\alpha}}\left\|\nabla^2V(x)\right\|.
\end{equation}
\end{definition}
It follows from definition that $U^\varepsilon \subseteq U^\varepsilon_{h,\alpha}$
and hence $L \leq L(h,\alpha)$. For $x\in U^\varepsilon$, Taylor's theorem yields
\begin{equation} \label{eq:V_taylor2}
	V(x + hY_h(x)) = V(x) + h \langle \nabla V(x), Y_h(x) \rangle + \tfrac{h^2}{2} Y_h(x)^\top \nabla^2 V(x + t hY_h(x)) Y_h(x)
\end{equation}
for some $t \in [0, 1]$. If we write $Y_h(x) = f(x) - \alpha \nabla V(x) + \mathcal{E}(x, h, \alpha)$ and substitute to~\eqref{eq:V_taylor2},
\begin{align} \label{eq:V_taylor3}
	V(x + hY_h(x)) &\leq V(x) - h\alpha \left|\nabla V(x)\right|^2 + h \left|\nabla V(x)\right| \left|\mathcal{E}(x, h, \alpha)\right| \notag \\
	&\qquad + \tfrac{h^2 L(h, \alpha)}{2} \left| f(x) - \alpha \nabla V(x) + \mathcal{E}(x, h, \alpha) \right|^2,
\end{align}
where $f(x) \perp \nabla V(x)$ is used. For a fixed step size $h$, a natural objective is to choose $\alpha$ so as to minimize the upper bound in~\eqref{eq:V_taylor3}. For a general one--step method, however, the remainder term $\mathcal{E}(x,h,\alpha)$ depends on both the chosen discretization method and the surrogate vector field. Hence, this optimization does not yield a method-agnostic closed-form rule for the gain. This reflects that gain selection is not intrinsic to the continuous feedback term alone, but is tied to the particular discrete realization of the surrogate dynamics. At the level of an unspecified one-step method, the local error structure is not fixed; once a discretization is specified, its remainder terms must be accounted for in the gain analysis. The Euler discretization is the canonical first-order explicit realization in which this remainder term vanishes, i.e., $\mathcal{E}(x,h,\alpha)=0$. Consequently, the gain dependence can be isolated explicitly as follows:
\begin{align} \label{eq:V_taylor4}
	V\bigl(x+h(f(x)-\alpha&\nabla V(x))\bigr) \notag \\ 
	& \leq V(x)-h\alpha |\nabla V(x)|^2 +\frac{h^2L(h,\alpha)}{2}|f(x)-\alpha\nabla V(x)|^2 \notag\\
	&= V(x)-h\alpha\left(1-\frac{h\alpha L(h,\alpha)}{2}\right)|\nabla V(x)|^2 +\frac{h^2L(h,\alpha)}{2}|f(x)|^2.
\end{align}

Now, let us introduce the scaled gain $\beta \coloneqq h\alpha$ so that the update rule in~\eqref{eq:feedback_discrete_system} can be written as
\begin{equation} \label{eq:feedback_discrete_system_optimal}
	x_{k+1}=x_k+hf(x_k)-\beta\nabla V(x_k).
\end{equation}
If $U^\varepsilon_{h,\alpha} = U^\varepsilon$, positive invariance of $U^\varepsilon$ under~\eqref{eq:feedback_discrete_system_optimal} is implied and $L(h,\alpha) = L$. In this case,~\eqref{eq:V_taylor4} reduces to
\begin{equation} \label{eq:V_taylor5}
	V\bigl(x+hf(x)-\beta \nabla V(x)\bigr) \leq V(x) + \left(-\beta+\frac{L}{2}\beta^2\right)|\nabla V(x)|^2 + \frac{h^2L}{2}|f(x)|^2,
\end{equation}
where the dissipative term on the right--hand side solely depends on $\beta$ but not on $h$. This leads to two questions: under which condition does $U^\varepsilon_{h,\alpha}$ coincide with $U^\varepsilon$, and, among admissible gains, which one minimizes the Taylor-based upper bound in~\eqref{eq:V_taylor5}. The next theorem answers both questions. Part~\ref{thm_optimal_a_item:1} identifies the admissible range $\beta\in\left(0,\frac{2}{L}\right)$, and part~\ref{thm_optimal_a_item:2} shows that the corresponding Taylor-based upper bound is minimized at $\beta=\frac{1}{L}$, equivalently $\alpha=\frac{1}{hL}$.
\begin{theorem} \label{thm:optimal_a}
	Suppose Feedback Integrator is implemented with Euler's method and fix $\varepsilon \in (0,\nu)$.
	\begin{enumerate}[label=(\roman*)]
		\item \label{thm_optimal_a_item:1} For every $\beta \in \left(0,\frac{2}{L}\right)$, there exists $h_\beta > 0$ such that, for every $0<h<h_\beta$, the scaled gain $\alpha=\frac{\beta}{h}$ satisfies $U^\varepsilon_{h,\alpha}\subseteq U^\varepsilon$, and therefore $U^\varepsilon_{h,\alpha} = U^\varepsilon$. Consequently, $U^\varepsilon$ is positively invariant under the discrete system~\eqref{eq:feedback_discrete_system} with $\alpha=\frac{\beta}{h}$.
		\item \label{thm_optimal_a_item:2} If $h>0$ and $\alpha>0$ satisfy $U^\varepsilon_{h,\alpha}\subseteq U^\varepsilon$ and if $\beta=h\alpha$, then for each $x\in U^\varepsilon$, the right--hand side of~\eqref{eq:V_taylor5} is minimized at $\beta = \frac{1}{L}$, or $\alpha=\frac{1}{hL}$. Moreover, by part~\ref{thm_optimal_a_item:1}, such gain choice is admissible for all sufficiently small $h$. This shows that $\alpha = \frac{1}{hL}$ is the optimal gain choice in the sense of minimizing the Taylor-based upper bound among all admissible gains.
	\end{enumerate}
\end{theorem}
\begin{proof}
	(i) Fix $\beta \in \left(0,\frac{2}{L}\right)$ and $x_k \in U^\varepsilon$. Define
	\begin{equation}
		x_{k+1} \coloneqq x_k + hf(x_k) - \beta \nabla V(x_k), \quad S_k \coloneqq \{\,x_k + t(x_{k+1}-x_k) \mid t \in [0,1]\,\}.
	\end{equation}
	Since $V$ is analytic, $U^\varepsilon$ is compact, and $U^\varepsilon$ contains no critical points of $V$ outside $V^{-1}(0)$, the {\L}ojasiewicz inequality yields $\mu>0$ and $\theta \in \left[\frac12,1\right)$ such that
	\begin{equation}
		|\nabla V(x)| \geq \mu V(x)^\theta, \qquad \forall x \in U^\varepsilon.
	\end{equation}
	For $\Updelta>0$, define $\overline{L} \coloneqq \sup\limits_{x\in U^{\varepsilon+\Updelta}}\left\|\nabla^2V(x)\right\|$ and $\overline{M}_f \coloneqq \sup\limits_{x\in U^{\varepsilon+\Updelta}}|f(x)|$. Then $\overline{L}\geq L$ and $\overline{L}\to L$ as $\Updelta\to0^+$. Since $\beta<2/L$, one can choose $\Updelta \in (0,\nu-\varepsilon)$ such that $\beta < 2/\overline{L}$. We first show that
	\begin{equation}
		S_k \subseteq U^{\varepsilon+\Updelta} \qquad\text{if}\qquad h < \sqrt{\frac{2\Updelta}{\overline{L}\,\overline{M}_f^2}}.
	\end{equation}
	To this end, let $\tau \coloneqq \inf\{\,t\ge0 \mid V(x_k+t(x_{k+1}-x_k)) \geq \varepsilon+\Updelta\,\}$ with the convention $\tau=\infty$ if the set is empty. For $t \leq \tau$, the point $x_k+t(x_{k+1}-x_k)$ lies in $U^{\varepsilon+\Updelta}$,
	so the same computation as in~\eqref{eq:V_taylor4} gives
	\begin{equation} \label{eq:V_segment1}
		V(x_k+t(x_{k+1}-x_k)) \leq V(x_k) - t\beta\left(1-\frac{t\beta\overline{L}}{2}\right)|\nabla V(x_k)|^2 + \frac{\overline{L}}{2}t^2h^2|f(x_k)|^2.
	\end{equation}
	Since $\beta<2/\overline{L}$, the second term on the right-hand side is nonpositive. Hence, 
	\begin{align} \label{eq:V_segment2}
		\sup\limits_{t\in[0,\min\{1,\tau\}]} V(x_k+t(x_{k+1}-x_k))
		&\leq V(x_k)+\frac{\overline{L}}{2}(\min\{1,\tau\})^2h^2|f(x_k)|^2 \notag\\
		&\leq V(x_k)+\frac{\overline{L}}{2}h^2|f(x_k)|^2.
	\end{align}
	If $\tau\le1$, then continuity gives $V(x_k+\tau(x_{k+1}-x_k))=\varepsilon+\Updelta$. But if $h < \sqrt{\frac{2\Updelta}{\overline{L}\,\overline{M}_f^2}}$, then the right-hand side of~\eqref{eq:V_segment2} is strictly smaller than $\varepsilon+\Updelta$, a contradiction. Thus, $\tau>1$, and therefore $S_k \subseteq U^{\varepsilon+\Updelta}$.

	Set
	\begin{equation}
		c_\beta^1 \coloneqq \beta\left(1-\frac{\beta\overline{L}}{2}\right)\mu^2,
		\qquad c_\beta^2 \coloneqq \frac{\overline{L}}{2}\overline{M}_f^2 .
	\end{equation}
	We now prove the stronger inclusion $S_k\subset U^\varepsilon$. Consider two cases.

	First, if $V(x_k)\leq \frac{\varepsilon}{2}$, then~\eqref{eq:V_segment1} and the nonpositivity of the dissipative term imply
	\begin{equation}
		V(x_k+t(x_{k+1}-x_k)) \leq \frac{\varepsilon}{2} + \frac{\overline{L}}{2}h^2\overline{M}_f^2 = \frac{\varepsilon}{2}+c_\beta^2h^2
	\end{equation}
	for all $t\in[0,1]$. Hence, $S_k\subseteq U^\varepsilon$ whenever $h \leq \sqrt{\frac{\varepsilon}{2c_\beta^2}}$.

	Assume next that $V(x_k)\in\left(\frac{\varepsilon}{2},\varepsilon\right]$. Since $t\in[0,1]$,
	\begin{equation}
		1-\frac{t\beta\overline{L}}{2}\geq 1-\frac{\beta\overline{L}}{2}.
	\end{equation}
	Using~\eqref{eq:V_segment1} and the {\L}ojasiewicz inequality,
	\begin{align}
		V(x_k+t(x_{k+1}-x_k))
		&\leq V(x_k) - t\beta\left(1-\frac{\beta\overline{L}}{2}\right)|\nabla V(x_k)|^2 + c_\beta^2 t^2 h^2 \notag \\
		&\leq \varepsilon - t c_\beta^1\left(\frac{\varepsilon}{2}\right)^{2\theta} + c_\beta^2 t^2 h^2.
	\end{align}
	Since $t^2\leq t$ on $[0,1]$,
	\begin{equation}
		V(x_k+t(x_{k+1}-x_k)) \leq \varepsilon - t\left[c_\beta^1\left(\frac{\varepsilon}{2}\right)^{2\theta} - c_\beta^2 h^2 \right].
	\end{equation}
	Therefore, $S_k\subseteq U^\varepsilon$ whenever $h \leq \sqrt{\frac{c_\beta^1}{c_\beta^2}}\left(\frac{\varepsilon}{2}\right)^\theta$.

	Combining the two cases, we conclude that $S_k \subseteq U^\varepsilon$ for every $x_k\in U^\varepsilon$, provided
	\begin{equation}
	h < h_\beta \coloneqq
	\min\left\{
		\sqrt{\frac{2\Updelta}{\overline{L}\,\overline{M}_f^2}},
		\sqrt{\frac{\varepsilon}{2c_\beta^2}},
		\sqrt{\frac{c_\beta^1}{c_\beta^2}}\left(\frac{\varepsilon}{2}\right)^\theta
	\right\}.
	\end{equation}
	Since $x_k\in U^\varepsilon$ was arbitrary, this proves that $U^\varepsilon_{h,\beta/h}\subseteq U^\varepsilon$ for all $0<h<h_\beta$. Taking $t=1$ implies $x_{k+1}\in U^\varepsilon$ whenever $x_k\in U^\varepsilon$, so $U^\varepsilon$ is positively invariant.

	(ii) Recall from definition that $U^\varepsilon\subseteq U^\varepsilon_{h,\alpha}$. Therefore, $U^\varepsilon_{h,\alpha}\subseteq U^\varepsilon$ implies $U^\varepsilon_{h,\alpha}=U^\varepsilon$, and hence $L(h,\alpha)=L$. For fixed $x$ and $h$, the $\beta$-dependent part in~\eqref{eq:V_taylor5} is the quadratic polynomial $q(\beta)\coloneqq -\beta+\frac{L}{2}\beta^2$, which is minimized at $q'(\beta)= -1+L\beta = 0 \iff \beta=\frac{1}{L}$. Equivalently, the minimizer is $\alpha=\frac{1}{hL}$. Since $\frac{1}{L}\in\left(0,\frac{2}{L}\right)$, part~\ref{thm_optimal_a_item:1} implies that this choice is admissible for all sufficiently small $h$. This proves the claimed optimality among admissible gains.
\end{proof}
\begin{remark}
	The Euler update in~\eqref{eq:feedback_discrete_system_optimal} can be viewed as a gradient descent with step size $\beta$, perturbed by an additive drift term $+hf(x_k)$. Indeed, by item~\ref{thm_optimal_a_item:1} of Theorem~\ref{thm:optimal_a}, for every $\beta\in\left(0,\frac{2}{L}\right)$ and all sufficiently small $h$, one has $U^\varepsilon_{h,\beta/h}\subseteq U^\varepsilon$. Hence $L(h,\beta/h)=L$, and~\eqref{eq:V_taylor4} yields, for $x_k\in U^\varepsilon$,
	\begin{equation}
		V(x_{k+1}) \leq V(x_k) -\beta\left(1-\frac{L}{2}\beta\right)\left|\nabla V(x_k)\right|^2 +\frac{h^2L}{2}\left|f(x_k)\right|^2.
	\end{equation}
	In the degenerate case $f\equiv0$, this reduces to
	\begin{equation}
		V(x_{k+1}) \leq V(x_k) -\beta\left(1-\frac{L}{2}\beta\right)\left|\nabla V(x_k)\right|^2,
	\end{equation}
	which leads to the standard convergence rate of gradient descent with step size $\beta$. Theorem~\ref{thm:optimal_a} further shows that even under the presence of the drift term $+hf(x_k)$, the iterates remain in $U^\varepsilon$.
\end{remark}
\begin{remark}
	For higher-order one-step methods, the remainder term $\mathcal{E}(x, h, \alpha)$ in~\eqref{eq:V_taylor3} depends on derivatives of the modified vector field $f-\alpha\nabla V$. Consequently, $\mathcal{E}(x, h, \alpha)$ generally contains mixed terms that are combinations of powers of $h$ and $\alpha$, which depend on both the integration method and the derivatives of $f$ and $V$. Hence, the upper bound in~\eqref{eq:V_taylor3} is no longer governed by a method- and dynamics- independent expression of gain, and the closed-form choice $\alpha=\frac{1}{hL}$ does not extend directly beyond Euler's method: higher--order contributions such as $O\left(h^2\alpha^2\right)$ remain $O\left(1\right)$ as $h \rightarrow 0$ and cannot be neglected even for small $h$. A comparable gain-selection rule for higher-order methods is left for future work.
\end{remark}
\subsection{Adaptive Gain Selection under Euler Discretization} \label{subsec:adaptive_gain}
Although Theorem~\ref{thm:optimal_a} identifies the fixed gain $\alpha=\frac{1}{hL}$ as the one that minimizes the Taylor-based upper bound, this choice has two practical limitations. First, evaluating
\begin{equation}
	L=\sup\limits_{x\in U^\varepsilon}\left\|\nabla^2V(x)\right\|
\end{equation}
is generally nontrivial, since it requires identifying the set $U^\varepsilon=V^{-1}([0,\varepsilon])$ and estimating the global maximum of $\|\nabla^2V\|$ on it. Second, when $\|\nabla^2V(x)\|$ varies substantially across $U^\varepsilon$, the global bound $L$ may be much larger than the local curvature scale along a given step. In such regions, the fixed gain $\alpha=\frac{1}{hL}$ becomes substantially smaller than the gain suggested by a local curvature bound, and therefore yields a conservative correction. This motivates a gain-selection rule that adapts to the local Hessian scale while retaining the positive-invariance property in Theorem~\ref{thm:optimal_a}.

For $x_k \in U^\varepsilon$, define
\begin{equation}
	B_{x_k} \coloneqq \max\left\{\left\|\nabla^2V(x_k)\right\|, H_{\min}\right\}, \quad \beta(x_k) \coloneqq \frac{1}{cB_{x_k}}, 
\end{equation}
where $c>1$ is a safety factor and $H_{\min}>0$ prevents the gain from becoming
arbitrarily large in regions where $\left\|\nabla^2V(x_k)\right\|$ is small. We define the Feedback Integrator under Euler's method with the adaptive gain $\beta(x_k)$ as
\begin{equation} \label{eq:feedback_discrete_system_adaptive}
	x_{k+1} \coloneqq x_k + h f(x_k) - \beta(x_k) \nabla V(x_k).
\end{equation}
To address the positive invariance of $U^\varepsilon$ under the update in~\eqref{eq:feedback_discrete_system_adaptive}, let us define
\begin{equation}
	K_{x_k} \coloneqq \left\{\,x \in U^\nu : \left\|\nabla^2V(x)\right\| \leq cB_{x_k}\,\right\},
\end{equation}
and
\begin{equation}
	U_{h,x_k,\beta} \coloneqq \bigl\{ x_k + u f(x_k) - v \nabla V(x_k) : u\in[0,h],\ v\in[0,\beta] \bigr\}.
\end{equation}
The set $U_{h,x_k,\beta(x_k)}$ contains the one-step segment joining $x_k$ and $x_{k+1}=x_k+h f(x_k)-\beta(x_k)\nabla V(x_k)$. Therefore, if $U_{h,x_k,\beta(x_k)}\subset K_{x_k}$, then $cB_{x_k}$ is a valid upper bound for $\left\|\nabla^2V\right\|$ on the region traversed by the next Euler step. Now recall that the proof of Theorem~\ref{thm:optimal_a} uses the Hessian only through an upper bound valid on the corresponding one-step segment.
\begin{lemma} \label{lem:local_one_step_bound}
	Fix $\varepsilon\in(0,\nu)$ and $x_k\in U^\varepsilon$. Let $\beta>0$ and $H>0$. Assume that
	\begin{equation}
		U_{h,x_k,\beta}\subset U^\nu, \qquad \sup\limits_{z\in U_{h,x_k,\beta}}\left\|\nabla^2V(z)\right\|\leq H.
	\end{equation}
	Then
	\begin{equation} \label{eq:local_one_step_bound}
		V\bigl(x_k+h f(x_k)-\beta \nabla V(x_k)\bigr) \leq V(x_k) -\beta\left(1-\frac{\beta H}{2}\right)|\nabla V(x_k)|^2 + \frac{H}{2}h^2|f(x_k)|^2.
	\end{equation}
	In particular, if $\beta<\frac{2}{H}$, then the dissipative term is nonpositive.
\end{lemma}
\begin{proof}
	Since $U_{h,x_k,\beta}$ contains the segment between $x_k$ and $x_k + hf(x_k)-\beta \nabla V(x_k)$, the conclusion follows from~\eqref{eq:V_taylor4} with $L(h,\beta/h)$ replaced by $H$.
\end{proof}
Lemma~\ref{lem:local_one_step_bound} reduces the analysis to two tasks. First, one must establish the inclusion
\begin{equation}
	U_{h,x_k,\beta(x_k)}\subset K_{x_k}
\end{equation}
for all $x_k\in U^\varepsilon$ and all sufficiently small $h$. Since $\beta(x_k)=\frac{1}{cB_{x_k}}$ is independent of $h$, the family $U_{h,x_k,\beta(x_k)}$ is monotone in $h$, and
the above inclusion is maintained as $h \rightarrow 0^+$. Under this inclusion, $cB_{x_k}$ is a valid Hessian bound on the one-step segment, and the choice $\beta(x_k)=\frac{1}{cB_{x_k}}$ satisfies $\beta(x_k)<\frac{2}{cB_{x_k}}$. Hence, the local one--step estimate in~\eqref{eq:local_one_step_bound} reduces to the same dissipative term as in the proof of Theorem~\ref{thm:optimal_a}. The remaining task is then to choose $h$ sufficiently small, uniformly in $x_k\in U^\varepsilon$, so that the $O(h^2)$ drift term in~\eqref{eq:local_one_step_bound} is also controlled. The following theorem addresses the two tasks and establishes the positive invariance of $U^\varepsilon$ under the update in~\eqref{eq:feedback_discrete_system_adaptive}.
\begin{theorem}\label{thm:adaptive_a}
	Assume that $\nabla^2V$ is $\Gamma$-Lipschitz on $U^\nu$. For any $c>1$, $H_{\min}>0$, and $\varepsilon\in(0,\nu)$, there exists $h_0>0$ such that, for every $0<h<h_0$, the sequence $\{x_k\}$ generated by the update~\eqref{eq:feedback_discrete_system_adaptive} from $x_0=x_I\in\Lambda$ satisfies
	\begin{equation}
		V(x_k)\leq \varepsilon \qquad \forall k \geq 0.
	\end{equation}
\end{theorem}
\begin{proof}
	Since $U^\varepsilon$ is a compact subset of the interior of $U^\nu$, there exists $r>0$ such that
	\begin{equation}
		\left(U^{\varepsilon}\right)_r \coloneqq \{x\in\mathbb R^n \mid \operatorname{dist}\left(x,U^{\varepsilon}\right)\leq r\}
		\subset U^\nu.
	\end{equation}
	Since $M$ in~\eqref{eq:gradient_bounds} is a class-$\mathcal{K}$ function, $M\left(\varepsilon'\right)\to 0$ as $\varepsilon'\to 0^+$. Hence, one can choose $\varepsilon'\in(0,\varepsilon)$ such that
	\begin{equation} \label{eq:adaptive_small_gradient}
		\frac{M\left(\varepsilon'\right)}{cH_{\min}}\leq \frac r2, \qquad
		\frac{\Gamma M\left(\varepsilon'\right)}{cH_{\min}} \leq \frac{c-1}{2}H_{\min}.
	\end{equation}
	Define 
	\begin{equation}
		F_{\varepsilon'} \coloneqq \sup\limits_{x\in U^{\varepsilon'}}|f(x)|, \qquad \overline B_{\varepsilon'} \coloneqq \sup\limits_{x\in U^{\varepsilon'}}B_x,
	\end{equation}
	and choose $h_{\mathrm{geo}}>0$ so that
	\begin{equation} \label{eq:adaptive_hgeo}
		h_{\mathrm{geo}}F_{\varepsilon'}\leq \frac r2, \qquad
		\Gamma h_{\mathrm{geo}}F_{\varepsilon'} \leq \frac{c-1}{2}H_{\min}.
	\end{equation}
	Fix any $x_k\in U^{\varepsilon'}$ and $h\in(0,h_{\mathrm{geo}}]$. For every
	$z\in U_{h,x_k,\beta(x_k)}$,
	\begin{equation} \label{eq:z_bound1}
		|z-x_k| \leq h|f(x_k)|+\beta(x_k)|\nabla V(x_k)| \leq hF_{\varepsilon'}+\frac{M\left(\varepsilon'\right)}{cH_{\min}}
		\leq r.
	\end{equation}
	Since $x_k\in U^{\varepsilon'}\subset U^{\varepsilon}$, it follows that
	\begin{equation}
		z\in \left(U^{\varepsilon}\right)_r \subset U^\nu \qquad \text{and} \qquad U_{h,x_k,\beta(x_k)} \subset U^\nu.
	\end{equation}
	Therefore, the $\Gamma$--Lipschitz continuity of $\nabla^2V$ gives
	\begin{align}
		\left\|\nabla^2V(z)\right\| &\leq \left\|\nabla^2V(x_k)\right\|+\Gamma |z-x_k| \notag \\
		&\leq \left\|\nabla^2V(x_k)\right\| +\Gamma hF_{\varepsilon'} +\frac{\Gamma M\left(\varepsilon'\right)}{cH_{\min}} \notag \\
		&\leq \left\|\nabla^2V(x_k)\right\|+(c-1)H_{\min}.
	\end{align}
	If $\left\|\nabla^2V(x_k)\right\|\geq H_{\min}$, then
	\begin{equation}
		\left\|\nabla^2V(z)\right\| \leq
		\left\|\nabla^2V(x_k)\right\|+(c-1)H_{\min} \leq c\left\|\nabla^2V(x_k)\right\| = cB_{x_k}.
	\end{equation}
	If $\left\|\nabla^2V(x_k)\right\|<H_{\min}$, then $B_{x_k}=H_{\min}$ and
	\begin{equation}
		\left\|\nabla^2V(z)\right\| \leq \left\|\nabla^2V(x_k)\right\|+(c-1)H_{\min} \leq cH_{\min} = cB_{x_k}.
	\end{equation}
	Hence,
	\begin{equation} \label{eq:adaptive_tube_inclusion}
		U_{h,x_k,\beta(x_k)}\subset K_{x_k} \qquad \forall x_k\in U^{\varepsilon'},\ \forall h\in(0,h_{\mathrm{geo}}].
	\end{equation}

	By Lemma~\ref{lem:local_one_step_bound} with
	\begin{equation} \label{eq:one_step_small_V_start}
		H=cB_{x_k}, \qquad \beta=\beta(x_k)=\frac{1}{cB_{x_k}},
	\end{equation}
	we obtain
	\begin{equation}
		V(x_{k+1}) \leq V(x_k) -\frac{1}{2cB_{x_k}}|\nabla V(x_k)|^2 +\frac{cB_{x_k}}{2}h^2|f(x_k)|^2.
	\end{equation}
	Since $B_{x_k}\leq \overline B_{\varepsilon'}$ and $|f(x_k)|\leq F_{\varepsilon'}$ on $U^{\varepsilon'}$,
	\begin{equation} \label{eq:adaptive_onestep_bound}
		V(x_{k+1}) \leq V(x_k) -\frac{1}{2c\overline B_{\varepsilon'}}|\nabla V(x_k)|^2 +\frac{c\overline B_{\varepsilon'}}{2}h^2F_{\varepsilon'}^2.
	\end{equation}
	Set
	\begin{equation}
		C_1 \coloneqq \frac{1}{2c\overline B_{\varepsilon'}}, \qquad
		C_2 \coloneqq \frac{c\overline B_{\varepsilon'}}{2}F_{\varepsilon'}^2.
	\end{equation}
	If $C_2=0$, then \eqref{eq:adaptive_onestep_bound} implies
	$V(x_{k+1})\leq V(x_k)$, and the conclusion follows immediately with $h_0 \coloneqq h_{\mathrm{geo}}$. Assume henceforth that $C_2>0$.

	Recall that the {\L}ojasiewicz inequality yields $\mu>0$ and $\theta \in \left[\frac12,1\right)$ such that
	\begin{equation} \label{eq:loja_adaptive}
		|\nabla V(x)| \geq \mu V(x)^\theta, \qquad \forall x \in U^{\varepsilon'}.
	\end{equation}

	We consider two cases. First, if $V(x_k)\leq \frac{\varepsilon'}{2}$, then
	\begin{equation}
		V(x_{k+1}) \leq \frac{\varepsilon'}{2}+C_2h^2.
	\end{equation}
	Hence, $V(x_{k+1})\leq \varepsilon'$ whenever
	\begin{equation}
		h\leq h_1 \coloneqq \sqrt{\frac{\varepsilon'}{2C_2}}.
	\end{equation}

	Assume next that $V(x_k)\in\left(\frac{\varepsilon'}{2},\varepsilon'\right]$. From~\eqref{eq:loja_adaptive},
	\begin{equation}
		|\nabla V(x_k)|\geq \mu V(x_k)^\theta \geq \mu\left(\frac{\varepsilon'}{2}\right)^\theta.
	\end{equation}
	Thus, \eqref{eq:adaptive_onestep_bound} yields
	\begin{equation}
		V(x_{k+1}) \leq \varepsilon' - C_1 \mu^2\left(\frac{\varepsilon'}{2}\right)^{2\theta} + C_2h^2.
	\end{equation}
	Therefore, $V(x_{k+1})\leq \varepsilon'$ whenever
	\begin{equation}
		h\leq h_2 \coloneqq \mu\sqrt{\frac{C_1}{C_2}} \left(\frac{\varepsilon'}{2}\right)^\theta.
	\end{equation}

	Finally, set
	\begin{equation}
		h_0 \coloneqq \min\{h_{\mathrm{geo}},h_1,h_2\}.
	\end{equation}
	Then
	\begin{equation} \label{eq:one_step_small_V_end}
		V(x_k)\leq \varepsilon' \quad\Longrightarrow\quad V(x_{k+1})\leq \varepsilon' \qquad \forall h\in(0,h_0).
	\end{equation}
	Since $x_0=x_I\in\Lambda$ and $V(x_0)=0$, induction gives
	\begin{equation}
		V(x_k)\leq \varepsilon' \qquad \forall k\geq 0.
	\end{equation}
	Because $\varepsilon'<\varepsilon$, the proof is complete.
\end{proof}
Theorem~\ref{thm:adaptive_a} treats a stepwise adaptive update, in which the gain is recomputed at every numerical step. As $h$ decreases, this increases the number of gain updates per unit physical time. The sequence $\left\|\nabla^2V(x_k)\right\|$ may be viewed as samples, at times $t_k=kh$, of the local Hessian scale $x \rightarrow \left\|\nabla^2V(x)\right\|$ along the computed trajectory. The variation of this Hessian scale occurs primarily along the trajectory in physical time, rather than being tied to the step index itself. Thus, reducing $h$ increases the sampling frequency of this variation, but does not, by itself, justify recomputing the gain at every numerical step. This motivates a time-periodic gain update.
\begin{algorithm}[t]
\caption{Feedback Integrator under Euler Discretization with Adaptive Gain} \label{alg:adaptive_gain_euler}
\begin{algorithmic}[1]
\Require Initial state $x_I \in \Lambda \subseteq U$, step size $h>0$, safety factor $c>1$, clip $H_{\min}>0$, max steps $N_{\max}$.
\Require Either: (i) stepwise update, or (ii) time-periodic update with period $T_{\mathrm{update}}>0$.
\Require $f(\cdot)$, $\nabla V(\cdot)$, $\mathrm{HessNorm}(x)=\left\|\nabla^2V(x)\right\|$.
\State $x_0 \gets x_I$, \quad $k \gets 0$
\If{time-periodic update is selected}
    \State $t_{\mathrm{next}} \gets 0$
\EndIf
\While{$k < N_{\max}$}
    \If{stepwise update is selected}
        \State $B \gets \max\{\mathrm{HessNorm}(x_k), H_{\min}\}$
        \State $\beta \gets \dfrac{1}{c\,B}$
    \ElsIf{$kh \geq t_{\mathrm{next}}$}
        \State $B \gets \max\{\mathrm{HessNorm}(x_k), H_{\min}\}$
        \State $\beta \gets \dfrac{1}{c\,B}$
        \State $t_{\mathrm{next}} \gets t_{\mathrm{next}} + T_{\mathrm{update}}$
    \EndIf
    \State $x_{k+1} \gets x_k + h f(x_k) - \beta \nabla V(x_k)$
    \State $k \gets k+1$
\EndWhile
\end{algorithmic}
\end{algorithm}
\begin{theorem} \label{thm:adaptive_a_period}
	Assume that $\nabla^2V$ is $\Gamma$--Lipschitz on $U^\nu$. Fix
	$\varepsilon\in(0,\nu)$, and assume that $V$ satisfies the {\L}ojasiewicz inequality
	\begin{equation}
		\left|\nabla V(x)\right|\geq \mu V(x)^\theta,\qquad \forall x\in U^\varepsilon,
	\end{equation}
	for some $\mu>0$ and $\theta\in\left[\frac{1}{2},1\right)$. Suppose moreover that the
	class--$\mathcal K$ function $M$ in~\eqref{eq:gradient_bounds} can be chosen so that
	\begin{equation}
		M(\delta)\leq C_M \delta^\theta,\qquad \forall \delta\in[0,\nu],
	\end{equation}
	for some constant $C_M>0$. Then, for any $c>1$ and $H_{\min}>0$, there exist $T_{\mathrm{update}}>0$ and $h_0>0$ such that, for every $0<h<h_0$, the sequence $\{x_k\}$ generated by Algorithm~\ref{alg:adaptive_gain_euler} with time-periodic gain update satisfies
	\begin{equation}
		V(x_k)\leq \varepsilon,\qquad \forall k\geq 0.
	\end{equation}
\end{theorem}
\begin{proof}
	Fix such a choice of $M$. Since $U^\varepsilon$ is a compact subset of the interior of $U^\nu$, there exists $r>0$ such that $\left(U^\varepsilon\right)_r\coloneqq\left\{x\in\mathbb{R}^n\mid \mathrm{dist}(x,U^\varepsilon)\leq r\right\}\subset U^\nu$.

	Define 
	\begin{equation}
		F_{\varepsilon} \coloneqq \sup\limits_{x\in U^{\varepsilon}}|f(x)|, \qquad \overline B_{\varepsilon} \coloneqq \sup\limits_{x\in U^{\varepsilon}}B_x,
	\end{equation}
	and set 
	\begin{equation}
		C_1\coloneqq\frac{1}{2c\overline{B}_\varepsilon}, \qquad C_2\coloneqq\frac{c\overline{B}_\varepsilon}{2}F_\varepsilon^2.
	\end{equation}
	Fix $\eta\in(0,\varepsilon]$, $x\in U^\eta$, and a gain
	\begin{equation}
		\beta_0\coloneqq\beta(x_{k_0})=\frac{1}{cB_{x_{k_0}}}
	\end{equation}
	for some $x_{k_0}\in U^\eta$. If
	\begin{equation}
		U_{h,x,\beta_0}\subset K_{x_{k_0}},
	\end{equation}
	then Lemma~\ref{lem:local_one_step_bound} with $H=cB_{x_{k_0}}\leq c\overline{B}_\varepsilon$ gives
	\begin{equation} \label{eq:thm_one_step_small_V}
		V\bigl(x+h f(x)-\beta_0\nabla V(x)\bigr) \leq V(x)-C_1|\nabla V(x)|^2 + C_2 h^2.
	\end{equation}
	If $C_2=0$, then the $O\left(h^2\right)$ term vanishes and~\eqref{eq:thm_one_step_small_V} yields
	\begin{equation}
		V(x)\leq \eta,\quad U_{h,x,\beta_0}\subset K_{x_{k_0}} \quad\Longrightarrow\quad V\bigl(x+h f(x)-\beta_0\nabla V(x)\bigr)\leq \eta.
	\end{equation}
	We now assume $C_2>0$ in what follows and define
	\begin{equation}
		\overline h_0(\eta) \coloneqq \min\left\{ \sqrt{\frac{\eta}{2C_2}}, \, \mu\sqrt{\frac{C_1}{C_2}}\left(\frac{\eta}{2}\right)^\theta \right\}.
	\end{equation}
	Through analogous steps as in~\eqref{eq:one_step_small_V_start}--\eqref{eq:one_step_small_V_end} in the proof of Theorem~\ref{thm:adaptive_a}, one has
	\begin{equation} \label{eq:one_step_eta_invariance}
		V(x)\leq \eta,\quad U_{h,x,\beta_0}\subset K_{x_{k_0}},\quad 0<h<\overline h_0(\eta) \quad\Longrightarrow\quad V\bigl(x+h f(x)-\beta_0\nabla V(x)\bigr)\leq \eta.
	\end{equation}

	Next, we prove, block by block, that
	\begin{equation}
		V(x_{k_0})\leq \eta \quad\Longrightarrow\quad V(x_{k_0+\ell})\leq \eta,\qquad \ell=0,1,\dots,n,
	\end{equation}
	where $n \coloneqq \left\lceil\frac{T_{\mathrm{update}}}{h}\right\rceil$ and $k_0$ is an index at which the gain is updated and the gain remains fixed on the block $\{k_0,\dots,k_0+n-1\}$. This is done by verifying the inclusion
	\begin{equation}
		U_{h,x_{k_0+\ell},\beta_0}\subset K_{x_{k_0}}
	\end{equation}
	uniformly along one gain-update block by choosing $T_{\mathrm{update}}$ and $h$ sufficiently small. This is the analogue of $h_{\mathrm{geo}}$ introduced in the proof of Theorem~\ref{thm:adaptive_a}.

	Since $\theta\in[1/2,1)$, by the definition of $\overline h_0$, there exist $\eta_0\in(0,\varepsilon)$ and $\overline C>0$ such that
	\begin{equation}
		\overline h_0(\eta)\geq \overline C\,\eta^\theta, \qquad \forall \eta\in(0,\eta_0).
	\end{equation}
	Choose $T_{\mathrm{update}}>0$ so small that
	\begin{equation} \label{eq:r_bound1}
		T_{\mathrm{update}} \left(F_\varepsilon+\frac{2C_M}{\overline C\,cH_{\min}}\right) \leq \frac{r}{2}, \qquad \Gamma T_{\mathrm{update}} \left(F_\varepsilon+\frac{2C_M}{\overline C\,cH_{\min}}\right) \leq \frac{c-1}{2}H_{\min}.
	\end{equation}
	Then choose $h_0\in\left(0,\frac{\overline C}{2}\eta_0^\theta\right)$ so small that
	\begin{equation} \label{eq:r_bound2}
		h_0
		\left(F_\varepsilon+\frac{2C_M}{\overline C\,cH_{\min}}\right) \leq \frac{r}{2}, \qquad \Gamma h_0 \left( F_\varepsilon+\frac{2C_M}{\overline C\,cH_{\min}} \right) \leq \frac{c-1}{2}H_{\min}.
	\end{equation}
	Fix $h\in(0,h_0)$ and let
	\begin{equation}
		\eta\coloneqq\left(\frac{2h}{\overline C}\right)^{1/\theta}.
	\end{equation}
	Then $\eta<\eta_0<\varepsilon$, and
	\begin{equation}
		h=\frac{\overline C}{2}\eta^\theta<\overline h_0(\eta).
	\end{equation}
	Moreover,
	\begin{equation}
		M(\eta)\leq C_M\eta^\theta = \frac{2C_M}{\overline C}h.
	\end{equation}

	Assume inductively that
	\begin{equation}
	V(x_{k_0+j})\leq \eta,\qquad j=0,1,\dots,\ell,
	\end{equation}
	for some $\ell\in\{0,\dots,n-1\}$. Let $z\in U_{h,x_{k_0+\ell},\beta_0}$. Since $\beta_0=\frac{1}{cB_{x_{k_0}}}\leq \frac{1}{cH_{\min}}$, as in~\eqref{eq:z_bound1}, we have
	\begin{align}
		|z-x_{k_0}|
		&\leq (\ell+1)\Bigl(hF_\varepsilon+\beta_0M(\eta)\Bigr) \notag \\
		&\leq n\left(hF_\varepsilon+\frac{M(\eta)}{cH_{\min}}\right) \notag \\
		&\leq \left(\frac{T_{\mathrm{update}}}{h}+1\right) \left(hF_\varepsilon+\frac{2C_M}{\overline C\,cH_{\min}}h\right) \notag \\
		&=(T_{\mathrm{update}}+h)\left(F_\varepsilon+\frac{2C_M}{\overline C\,cH_{\min}}\right) \notag \\
		&\leq r,
	\end{align}
	where the last inequality follows from~\eqref{eq:r_bound1} and~\eqref{eq:r_bound2}. Hence,
	\begin{equation}
	z\in \left(U^\varepsilon\right)_r \subset U^\nu.
	\end{equation}
	Using the $\Gamma$--Lipschitz continuity of $\nabla^2V$ on $U^\nu$,
	\begin{align}
		\left\|\nabla^2V(z)\right\|
		&\leq \left\|\nabla^2V(x_{k_0})\right\| + \Gamma |z-x_{k_0}| \notag \\
		&\leq \left\|\nabla^2V(x_{k_0})\right\| + \Gamma (T_{\mathrm{update}}+h) \left(F_\varepsilon+\frac{2C_M}{\overline C\,cH_{\min}}\right) \notag \\
		&\leq \left\|\nabla^2V(x_{k_0})\right\| + (c-1)H_{\min}
		\leq cB_{x_{k_0}}.
	\end{align}
	Therefore,
	\begin{equation}
		U_{h,x_{k_0+\ell},\beta_0}\subset K_{x_{k_0}}.
	\end{equation}
	Since $h<\overline h_0(\eta)$, the one--step estimate established above implies
	\begin{equation}
		V(x_{k_0+\ell+1})\leq \eta.
	\end{equation}
	This proves the block claim.

	Finally, since $x_0=x_I\in\Lambda$, we have $V(x_0)=0<\eta$. Applying the block
	claim inductively over all gain-update blocks yields
	\begin{equation}
		V(x_k)\leq \eta<\varepsilon,\qquad \forall k\geq 0.
	\end{equation}
	This completes the proof.
\end{proof}
Over a fixed interval $[0,T]$, the number of time-periodic gain updates is at most $\left\lceil \frac{T}{T_{\mathrm{update}}}\right\rceil+1$, and hence does not increase as $h \rightarrow 0^+$. Consequently, the additional overhead from recomputing the gain, measured relative to the total integration cost, vanishes as $h \rightarrow 0^+$; in this sense, its computational cost approaches that of the fixed-gain feedback integrator in the small-step regime. At the same time, the method avoids an a priori estimate of the global Hessian bound $L$. To make Theorem~\ref{thm:adaptive_a_period} directly applicable, we now verify its hypotheses in the standard sum-of-squares constant-rank setting.
\begin{proposition}[The sum-of-squares constant-rank case implies $\theta=\frac12$] \label{prop:theta_q_half}
	Under Assumption~\ref{assumption:01}, suppose in addition that
	\begin{equation}\label{eq:V_sos_form}
		V(x)=\frac{1}{2} |g(x)|^2, \qquad \Lambda = V^{-1}(0)=g^{-1}(0),
	\end{equation}
	for a smooth map $g:U\to\mathbb{R}^r$, and that $Dg$ has constant rank in a neighborhood of $\Lambda$.
	Then the following hold.

	\begin{enumerate}[label=(\roman*)]
	\item \label{prop:theta_q_half_item1}
	$\Lambda$ is a Morse--Bott critical submanifold of $V$, and $V$ is a Morse--Bott function along $\Lambda$. More precisely, for every $x\in\Lambda$,
	\begin{equation}
		\nabla V(x)=0,\qquad
		\nabla^2V(x)=Dg(x)^\top Dg(x),\qquad
		\ker\bigl(\nabla^2V(x)\bigr)=T_x\Lambda.
	\end{equation}

	\item \label{prop:theta_q_half_item2}
	There exist a neighborhood $\mathcal N$ of $\Lambda$ and constants $c_1,c_2,c_3,c_4>0$ such that, for all $x\in\mathcal N$,
	\begin{equation}
		c_1\,\operatorname{dist}(x,\Lambda)^2 \leq V(x)\leq c_2\,\operatorname{dist}(x,\Lambda)^2,
	\end{equation}
	and
	\begin{equation}
		c_3\,\operatorname{dist}(x,\Lambda)\leq |\nabla V(x)|\leq c_4\,\operatorname{dist}(x,\Lambda).
	\end{equation}

	\item \label{prop:theta_q_half_item3}
	There exist $C_M>0$ and a class--$\mathcal K$ function $M$ in~\eqref{eq:gradient_bounds} such that
	\begin{equation}
		M(\delta)\leq C_M\sqrt{\delta},\qquad \forall \delta\in[0,\nu].
	\end{equation}

	\item \label{prop:theta_q_half_item4}
	For every $\varepsilon\in(0,\nu)$, there exists $\mu_\varepsilon>0$ such that
	\begin{equation}
		|\nabla V(x)|\geq \mu_\varepsilon \sqrt{V(x)},\qquad \forall x\in U^\varepsilon.
	\end{equation}
	\end{enumerate}

	Consequently, the hypothesis on $M$ in Theorem~\ref{thm:adaptive_a_period} and the {\L}ojasiewicz inequality there hold with $\theta=\frac12$.
\end{proposition}
\begin{proof}
	Since
	\begin{equation}
		V(x)=\frac12\sum_{i=1}^r g_i(x)^2,
	\end{equation}
	one has
	\begin{equation}
		\nabla V(x)=\sum_{i=1}^r g_i(x)\nabla g_i(x)=Dg(x)^\top g(x),
	\end{equation}
	and
	\begin{equation}\label{eq:full_hessian_formula}
	\nabla^2V(x)=Dg(x)^\top Dg(x)+\sum_{i=1}^r g_i(x)\nabla^2 g_i(x).
	\end{equation}
	Hence, for every $x\in\Lambda=g^{-1}(0)$,
	\begin{equation}
		\nabla V(x)=0, \qquad \nabla^2V(x)=Dg(x)^\top Dg(x).
	\end{equation}

	Because $Dg$ has constant rank in a neighborhood of $\Lambda$, the constant-rank theorem
	implies that $\Lambda=g^{-1}(0)$ is a smooth submanifold. Moreover, for $x\in\Lambda$,
	\begin{equation}
		T_x\Lambda=\ker Dg(x).
	\end{equation}
	Using $\nabla^2V(x)=Dg(x)^\top Dg(x)$, we obtain
	\begin{equation}
		\ker\bigl(\nabla^2V(x)\bigr)
		=
		\ker\bigl(Dg(x)^\top Dg(x)\bigr)
		=
		\ker Dg(x)
		=
		T_x\Lambda.
	\end{equation}
	Since $Dg(x)^\top Dg(x)\succeq 0$, its restriction to the normal space $N_x\Lambda$ is positive definite. Therefore, $\Lambda$ is a smooth critical submanifold of $V$, and $V$ is a Morse--Bott function along $\Lambda$. This proves~\ref{prop:theta_q_half_item1}.

	We now prove~\ref{prop:theta_q_half_item2}. By the Morse--Bott lemma (with index $0$, since $V\geq 0$ and $V|_\Lambda=0$), for each $p\in\Lambda$ there exist local coordinates $(u,v)\in\mathbb{R}^d\times\mathbb{R}^{n-d}$ where $d=\dim\Lambda$, centered at $p$, such that
	\begin{equation}
		\Lambda=\{v=0\},\qquad V(u,v)=|v|^2
	\end{equation}
	in those coordinates; see, for example, Theorem~2 in~\cite{banyaga2004proof}. Since $\Lambda=V^{-1}(0)$ is a closed subset of the compact set $V^{-1}([0,\nu])$, it is compact. By the tubular neighborhood theorem~\cite{lee2018riemannian}, after shrinking the neighborhood if necessary, there exist a neighborhood $\mathcal{N}$ of $\Lambda$, a radius $\rho>0$, a smooth nearest-point projection
	\begin{equation}
		\Pi:\mathcal N\to\Lambda,
	\end{equation}
	and a vector field $\xi: \mathcal N \rightarrow \mathbb R^n$ such that every $x\in\mathcal N$ can be written uniquely as
	\begin{equation} \label{eq:x:as:sum}
		x=p+\xi_x,\qquad p=\Pi(x)\in\Lambda,\qquad \xi_x \in N_p\Lambda \subset \mathbb{R}^n,\qquad |\xi_x|<\rho,
	\end{equation}
	and
	\begin{equation} \label{eq:dist:x:to:Lambda}
		\operatorname{dist}(x,\Lambda)=|\xi_x|.
	\end{equation}
	Let
	\begin{equation}
		H_p \coloneqq \nabla^2V(p)|_{N_p\Lambda}
	\end{equation}
	for $p$ in $\Lambda$. Since $V$ is Morse--Bott along $\Lambda$ and $V\ge0$ with $V|_\Lambda=0$, the quadratic form $H_p$ is positive definite on $N_p\Lambda$. By compactness of $\Lambda$, there exist constants $\lambda_-,\lambda_+>0$ such that
	\begin{equation}
		\lambda_-|\zeta|^2 \leq \left\langle H_p\zeta, \zeta \right\rangle \leq \lambda_+|\zeta|^2
	\end{equation}
	for all $p\in\Lambda$ and all $\zeta\in N_p\Lambda$.
		
	Shrinking $\rho>0$ if necessary, continuity of $\nabla^2V$ gives
	\begin{equation}
		\left\langle \nabla^2V(p+s\xi)\xi, \xi \right\rangle \geq \frac{\lambda_-}{2}|\xi|^2
	\end{equation}
	for all $p\in\Lambda$, $s\in[0,1]$, and $\xi\in N_p\Lambda$ with $|\xi|<\rho$. Since $V(p)=0$ and $\nabla V(p)=0$, Taylor's formula yields
	\begin{equation}
		V(p+\xi) = \int_0^1(1-s) \left\langle \nabla^2V(p+s\xi)\xi, \xi \right\rangle\,ds.
	\end{equation}
	Hence,
	\begin{equation} \label{eq:V:geq}
		V(p+\xi)\geq \frac{\lambda_-}{4}|\xi|^2.
	\end{equation}
	Similarly, since $\nabla^2V$ is bounded on $\mathcal N$, there exists $C_V>0$ such that
	\begin{equation} \label{eq:V:leq}
		V(p+\xi)\leq C_V|\xi|^2.
	\end{equation}
	By~\eqref{eq:x:as:sum},~\eqref{eq:dist:x:to:Lambda},~\eqref{eq:V:geq}, and~\eqref{eq:V:leq}, we obtain
	\begin{equation}
		c_1\operatorname{dist}(x,\Lambda)^2 \leq V(x) \leq c_2\operatorname{dist}(x,\Lambda)^2
	\end{equation}
	for some $c_1,c_2>0$ and all $x\in\mathcal N$.

	It remains to estimate the gradient. Since $\nabla V(p)=0$,
	\begin{equation}
		\nabla V(p+\xi) = \int_0^1 \nabla^2V(p+s\xi)\xi\,ds.
	\end{equation}
	The upper bound
	\begin{equation} \label{eq:nabla:V:leq}
		|\nabla V(p+\xi)|\leq C|\xi|
	\end{equation}
	follows from boundedness of $\nabla^2V$ on $\mathcal N$. For the lower bound, let $P_p$ denote the orthogonal projection onto $N_p\Lambda$. Since $\nabla^2V(p)\xi=H_p\xi\in N_p\Lambda$, shrinking $\mathcal N$ if necessary gives
	\begin{equation}
		\left| P_p\nabla V(p+\xi)-H_p\xi \right| \leq \frac{\lambda_-}{2}|\xi|.
	\end{equation}
	Hence,
	\begin{equation} \label{eq:nabla:V:geq}
		|\nabla V(p+\xi)| \geq |P_p\nabla V(p+\xi)| \geq |H_p\xi|-\frac{\lambda_-}{2}|\xi| \geq \frac{\lambda_-}{2}|\xi|.
	\end{equation}
	By~\eqref{eq:x:as:sum},~\eqref{eq:dist:x:to:Lambda},~\eqref{eq:nabla:V:leq}, and~\eqref{eq:nabla:V:geq},
	\begin{equation}
		c_3\operatorname{dist}(x,\Lambda) \leq |\nabla V(x)| \leq c_4\operatorname{dist}(x,\Lambda)
	\end{equation}
	for some $c_3,c_4>0$ and all $x\in\mathcal N$. This proves~\ref{prop:theta_q_half_item2}.
	
	Define
	\begin{equation}
		S(\delta)\coloneqq \sup\limits_{x \in U^\delta}|\nabla V(x)|, \qquad \delta\in[0,\nu].
	\end{equation}
	Choose $\delta_0\in(0,\nu)$ so that $U^{\delta_0}\subset \mathcal N$. If $0\leq \delta\leq \delta_0$ and $V(x)\leq \delta$, then $x\in \mathcal N$, so by
	\ref{prop:theta_q_half_item2},
	\begin{equation}
		|\nabla V(x)|
		\leq c_4\,\operatorname{dist}(x,\Lambda)
		\leq \frac{c_4}{\sqrt{c_1}}\,\sqrt{V(x)}
		\leq \frac{c_4}{\sqrt{c_1}}\,\sqrt{\delta}.
	\end{equation}
	Hence,
	\begin{equation}
		S(\delta)\leq \frac{c_4}{\sqrt{c_1}}\,\sqrt{\delta},
	\qquad 0\leq \delta\leq \delta_0.
	\end{equation}
	For $\delta\in[\delta_0,\nu]$, compactness of $U^\nu$ implies $S(\nu) < \infty$, so
	\begin{equation}
		S(\delta)\leq S(\nu)\leq \frac{S(\nu)}{\sqrt{\delta_0}}\,\sqrt{\delta}.
	\end{equation}
	Therefore, with
	\begin{equation}
		C_M \coloneqq \max\left\{\frac{c_4}{\sqrt{c_1}}, \frac{S(\nu)}{\sqrt{\delta_0}} \right\},
	\end{equation}
	we have
	\begin{equation}
		S(\delta)\leq C_M\sqrt{\delta},\qquad \forall \delta\in[0,\nu].
	\end{equation}
	Hence, the class--$\mathcal K$ function
	\begin{equation}
		M(\delta)\coloneqq C_M\sqrt{\delta}
	\end{equation}
	satisfies~\eqref{eq:gradient_bounds}. This proves
	item~\ref{prop:theta_q_half_item3}.

	Now fix $\varepsilon\in(0,\nu)$. On $U^\varepsilon\cap \mathcal N$, item
	\ref{prop:theta_q_half_item2} gives
	\begin{equation} \label{eq:thm_theta_q_half_item4_proof1}
		|\nabla V(x)| \geq c_3\,\operatorname{dist}(x,\Lambda) \geq \frac{c_3}{\sqrt{c_2}}\,\sqrt{V(x)}.
	\end{equation}
	Consider the compact set $U^\varepsilon\setminus \mathcal N$. If $U^\varepsilon\setminus \mathcal N=\emptyset$, then~\eqref{eq:thm_theta_q_half_item4_proof1} proves~\ref{prop:theta_q_half_item4}. Otherwise, $(U^\varepsilon\setminus \mathcal N) \cap \Lambda=\emptyset$, and Assumption~\ref{assumption:01},~\ref{assumption:01:A3} implies that $|\nabla V|$ has no zeros on $U^\varepsilon\setminus \mathcal N$. Hence,
	\begin{equation}
	\gamma_\varepsilon\coloneqq\min_{x\in U^\varepsilon\setminus \mathcal N} |\nabla V(x)|>0.
	\end{equation}
	Since $V(x)\leq \varepsilon$ on $U^\varepsilon$, it follows that
	\begin{equation}
	|\nabla V(x)|
	\geq \gamma_\varepsilon
	\geq \frac{\gamma_\varepsilon}{\sqrt{\varepsilon}}\,\sqrt{V(x)},
	\qquad x\in U^\varepsilon\setminus \mathcal N.
	\end{equation}
	Therefore, with
	\begin{equation}
	\mu_\varepsilon
	\coloneqq
	\min\left\{
	\frac{c_3}{\sqrt{c_2}},
	\frac{\gamma_\varepsilon}{\sqrt{\varepsilon}}
	\right\},
	\end{equation}
	we obtain
	\begin{equation}
	|\nabla V(x)|\geq \mu_\varepsilon \sqrt{V(x)},
	\qquad \forall x\in U^\varepsilon.
	\end{equation}
	This proves~\ref{prop:theta_q_half_item4}. The final claim follows immediately
	from items~\ref{prop:theta_q_half_item3} and~\ref{prop:theta_q_half_item4}.
\end{proof}
%
%
\section{Numerical Demonstration} \label{sec:demo}
We present numerical demonstration results of Feedback Integrator conducted on (i) free rigid body motion on $\operatorname{SO}(3)$, (ii) the Kepler problem, and (iii) a perturbed Kepler problem with rotational symmetry. The conditions in Assumption~\ref{assumption:01} of all Lyapunov functions introduced throughout the demonstrations are verified in~\cite{chang2016feedback}. Comparisons are made among feedback integrators with Euler discretization as the baseline integration scheme, with gains $\alpha = 1$ (unity), $\alpha = \frac{1}{hL}$ (Theorem~\ref{thm:optimal_a}), adaptive gain under time-periodic update with Algorithm~\ref{alg:adaptive_gain_euler}, and with standard benchmark methods in~\cite{Hairer2006}.

For the Taylor-based fixed gain rule $\alpha = \frac{1}{hL}$ suggested by Theorem~\ref{thm:optimal_a}, we use a trajectory-based estimate of the Hessian scale, obtained by recording the maximum observed $\left\|\nabla^2 V(x_k)\right\|$ along a unity-gain feedback trajectory over one representative period. This estimate is used only for numerical comparison. Strictly speaking, efficient and accurate estimation of Lipschitz constant is another problem per se, especially when there is no additional knowledge on the dynamics such as periodicity. Not requiring such estimation step is one of the biggest benefit of adaptive gain selection as outlined in Section~\ref{subsec:adaptive_gain}.

Throughout all demonstrations with adaptive gain, we set $c = 1.1$ and $H_{\min} = 10^{-10}$, and Frobenius norm is used instead of matrix 2--norm for gain calculation. A run is marked as divergent when $V(x_k)$ exceeds $10^5$ or the numerical trajectory leaves the displayed bounded region. All simulations are conducted on MacBook Air M2 with C++, and codes are available at: \url{https://github.com/johnbae1901/Feedback-Integrator}.
\subsection{Free Rigid Body Motion in $\operatorname{SO}(3)$}
We consider the free rigid body dynamics as follows, 
\begin{subequations}\label{eq:dynamics_rigid}
\begin{equation} 
	\dot{R} = R \hat{\Omega}
\end{equation}
\begin{equation} 
	\dot{\Omega} = \mathbb{I}^{-1}((\mathbb{I}\Omega) \times \Omega)
\end{equation}
\end{subequations}
where $(R, \Omega) \in \operatorname{SO}(3) \times \mathbb{R}^3$, $\mathbb{I}$ is the moment of inertia matrix, and $\hat{\Omega}\in\mathfrak{so}(3)$ denotes the skew-symmetric matrix satisfying $\hat{\Omega}a=\Omega\times a$ for all $a\in\mathbb R^3$. To apply feedback integrator, we assume that the system is defined through the same expressions in~\eqref{eq:dynamics_rigid} in $\mathbb{R}^{3\times 3} \times \mathbb{R}^3$. The first integrals of this system are kinetic energy and spatial angular momentum, represented as follows respectively.
\begin{equation}
	E(\Omega) = \frac{1}{2}\Omega^\top \mathbb{I} \Omega, \qquad \pi(R, \Omega) = R\mathbb{I}\Omega
\end{equation}
For initial values $R_I \in \operatorname{SO}(3)$ and $\Omega_I \in \mathbb{R}^3 \setminus \{(0, 0, 0)\}$, let us define $E_I \coloneqq E(\Omega_I)$ and $\pi_I \coloneqq \pi(R_I, \Omega_I)$. Define an open set $U\coloneqq \left\{ (R, \Omega) \in \mathbb{R}^{3\times 3} \times \mathbb{R}^3 \mid \det(R) > 0 \right\}$ and the Lyapunov function $V: U \rightarrow \mathbb{R}_{\geq 0}$ as 
\begin{equation}
	V(R, \Omega) \coloneqq \frac{k_0}{4}\left\| R^\top R - I \right\|_F^2 + \frac{k_1}{2}\left| E(\Omega) - E_I \right|^2 + \frac{k_2}{2}	\left| \pi(R, \Omega) - \pi_I \right|^2,
\end{equation}
for constants $k_0, k_1, k_2 > 0$ and $3\times 3$ identity matrix $I$.

Throughout the simulations, we use $\mathbb{I} = \mathrm{diag}(3, 2, 1)$, $R_I = I$, and $\Omega_I = (1, 1, 1)$, which correspond to $E_I = 3$ and $\pi_I = (3, 2, 1)$. For the Lyapunov function, we use $k_0 = 50$, $k_1 = 100$, $k_2 = 50$. The Lipschitz constant of $\nabla V$ is estimated as $L \approx 1986.0$. For feedback integrator with adaptive gain selection, the gain update is set to be done every $T_{\mathrm{update}} = 30$ seconds. We present the results for $h \in \left[10^{-7}, 10^{-1}\right]$, and comparisons are made with a benchmark framework, the \emph{Strang Splitting method}~\cite{Hairer2006}.

\begin{figure}[t]
  \centering
  \begin{subfigure}[t]{0.96\linewidth}
    \includegraphics[width=\linewidth]{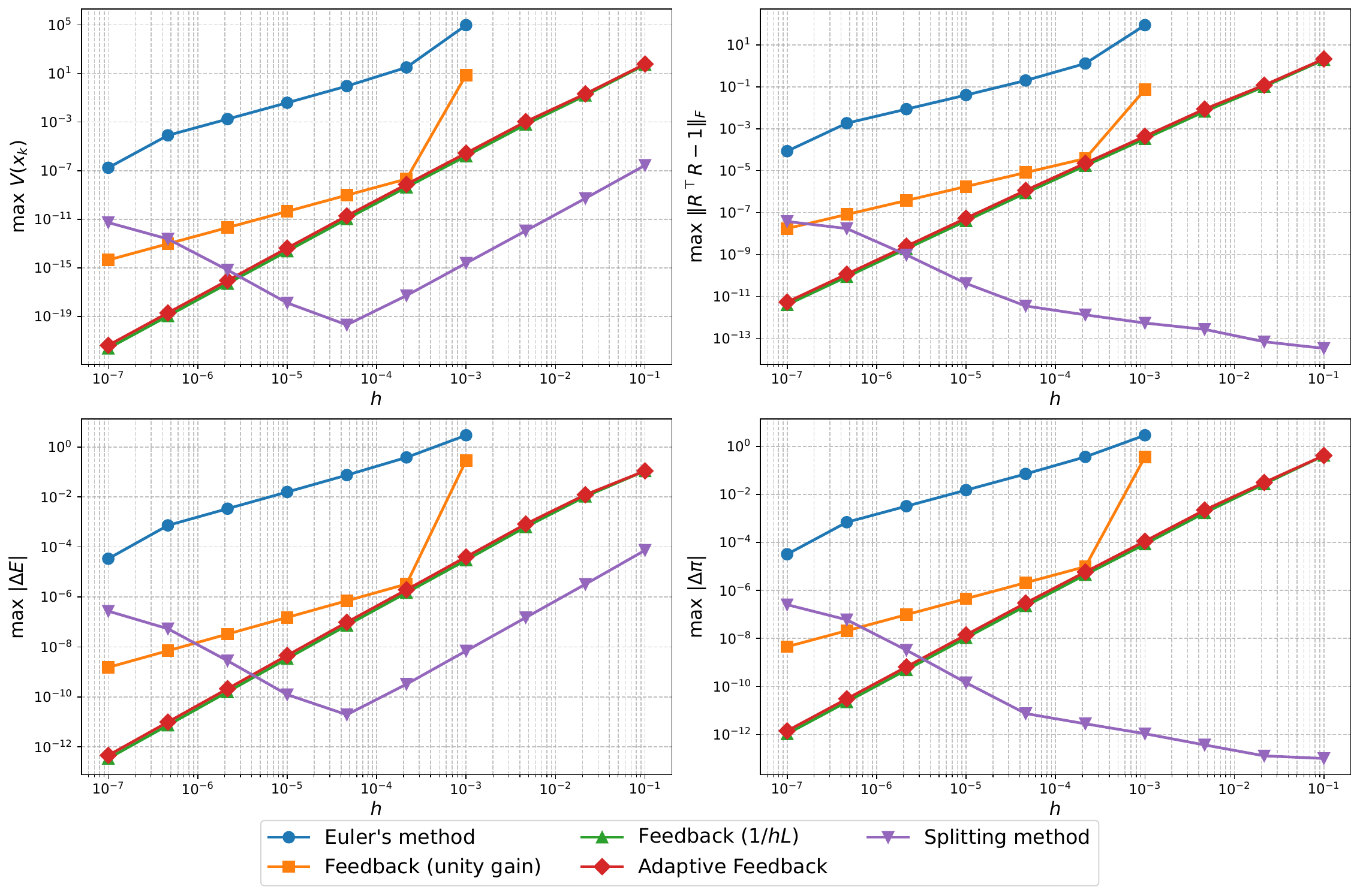}
    \caption{Accuracy vs.\ $h$}
  \end{subfigure}\hfill
  \begin{subfigure}[t]{0.6\linewidth}
    \includegraphics[width=\linewidth]{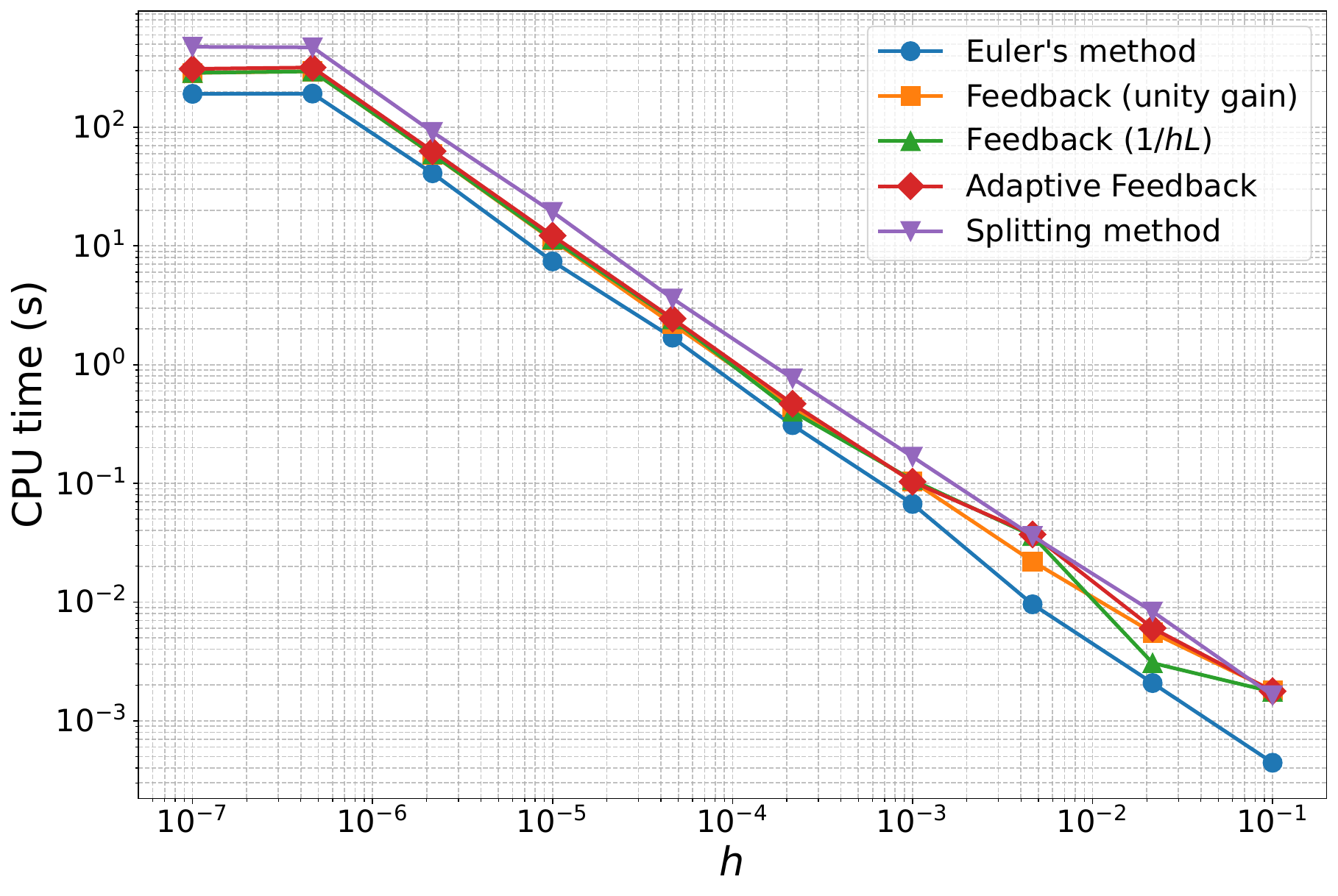}
    \caption{CPU time vs.\ $h$}
  \end{subfigure}
  \caption[Accuracy and cost vs.\ step size $h$]%
  {\textbf{Accuracy results of free rigid body motion in $\operatorname{SO}(3)$.} Integration with Euler's method on $[0, 1000]$. (a) maximum $V(x_k)$ and deviation of first integrals along the trajectories. (b) CPU time dedicated for each integration scheme.} \label{fig:rigid_body_results}
\end{figure}

\begin{figure}[t]
  \centering
  \begin{subfigure}[t]{0.32\linewidth}
    \includegraphics[width=\linewidth]{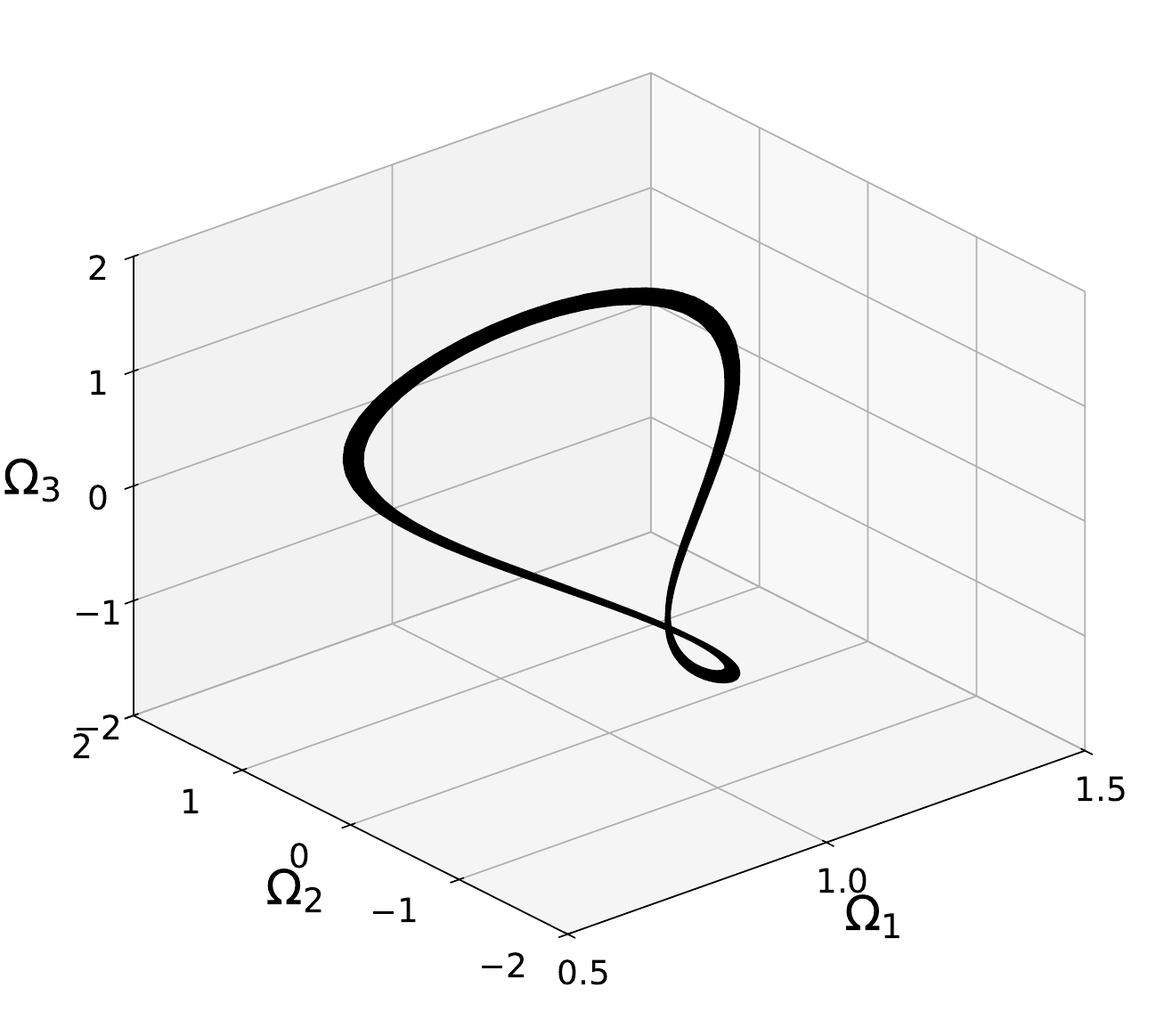}
    \caption{Euler's method}
  \end{subfigure}\hfill
  \begin{subfigure}[t]{0.32\linewidth}
    \includegraphics[width=\linewidth]{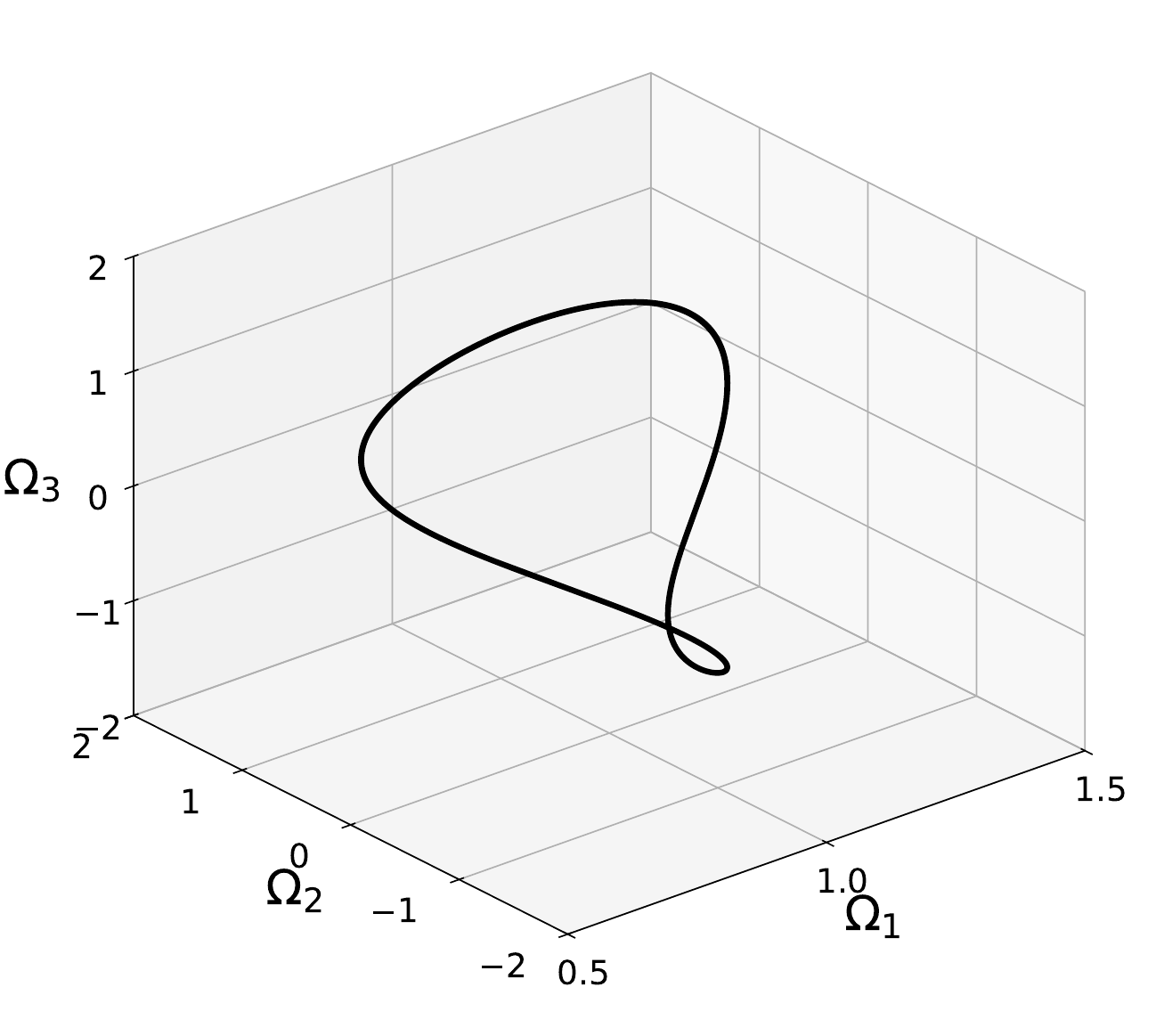}
    \caption{Feedback ($\alpha = 1$)}
  \end{subfigure}\hfill
  \begin{subfigure}[t]{0.32\linewidth}
    \includegraphics[width=\linewidth]{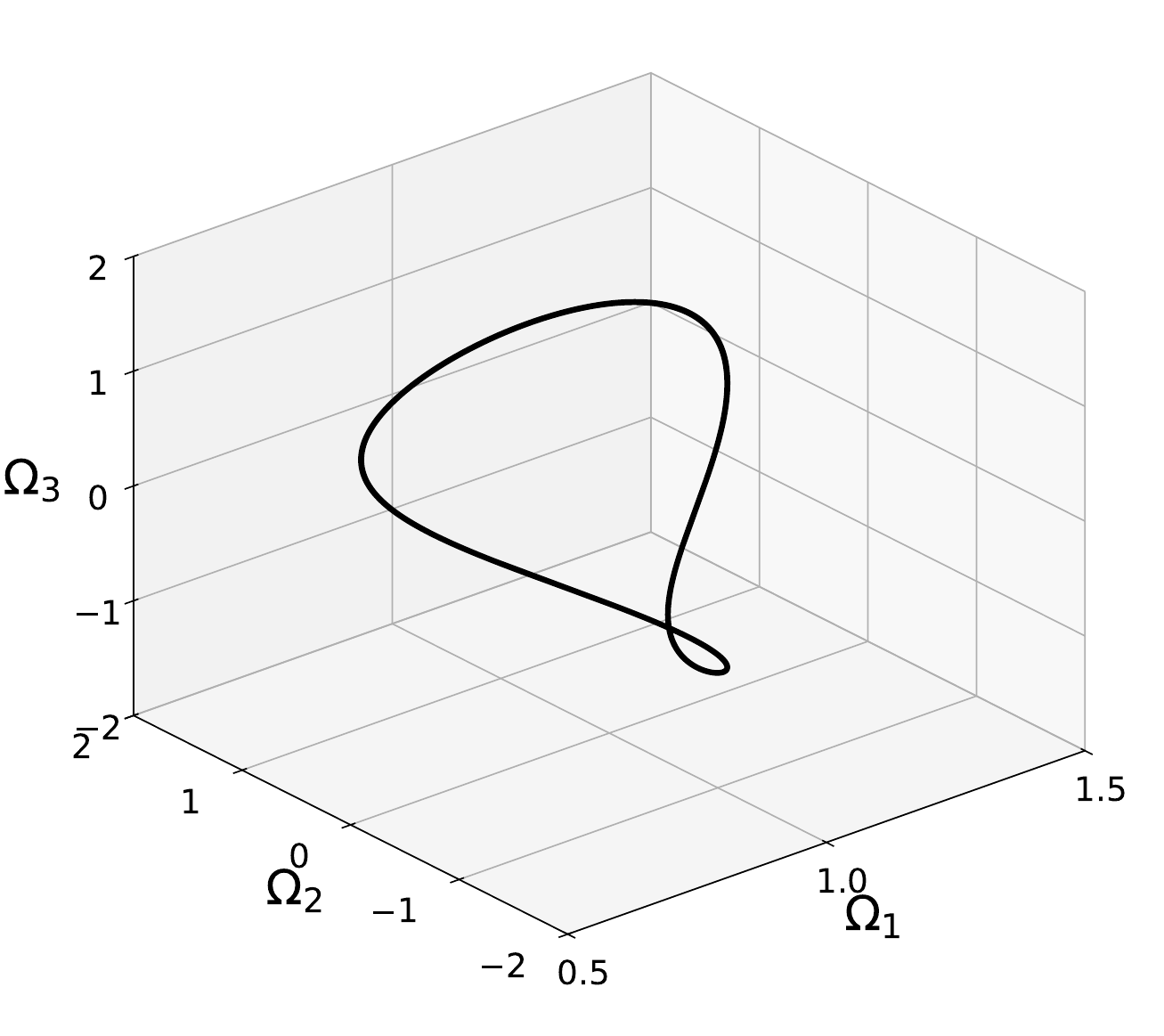}
    \caption{Feedback ($\alpha = 1/hL$)}
  \end{subfigure}\hfill
  \begin{subfigure}[t]{0.48\linewidth}
    \includegraphics[width=\linewidth]{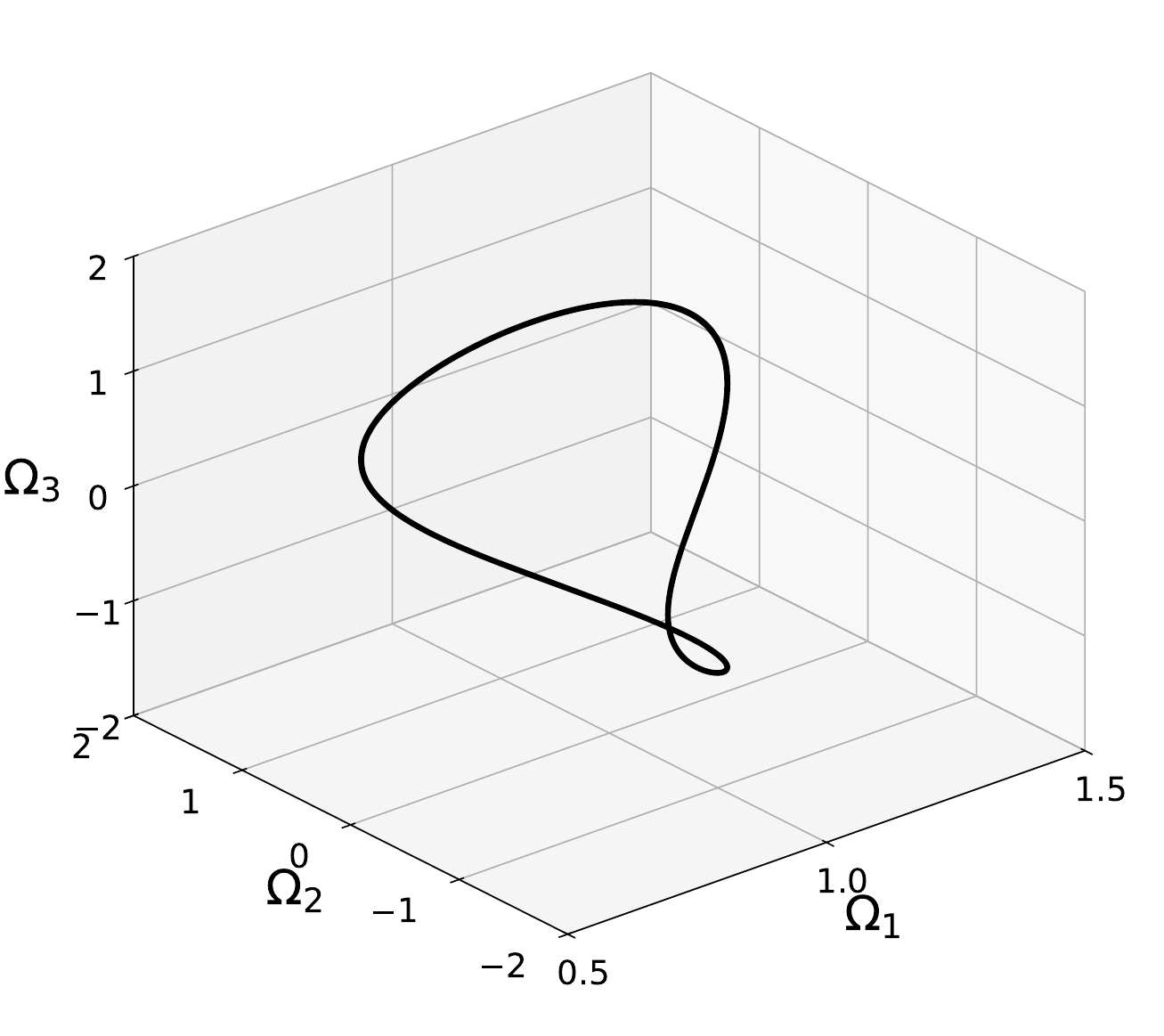}
    \caption{Adaptive feedback}
  \end{subfigure}\hfill
  \begin{subfigure}[t]{0.48\linewidth}
    \includegraphics[width=\linewidth]{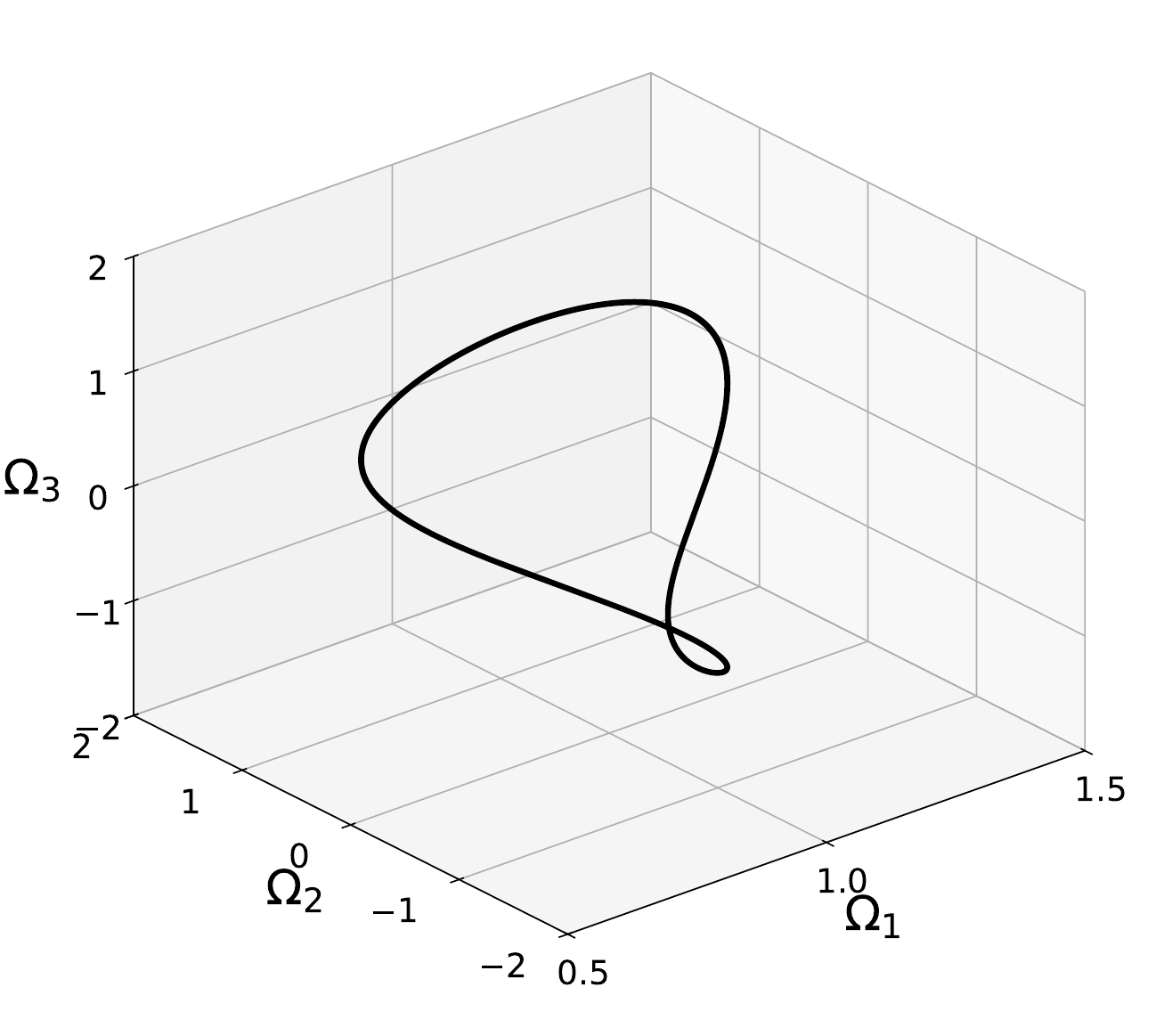}
    \caption{Splitting method}
  \end{subfigure}\hfill
  \caption[Rigid Body Trajectories]
  {\textbf{Trajectories of the body angular velocities of free rigid body motion in $\operatorname{SO}(3)$.} Integration on $[0, 1000]$ with $h = 10^{-4}$. Feedback integrators are implemented with Euler's method.} \label{fig:rigid_body_traj}
\end{figure}

Figure~\ref{fig:rigid_body_results} summarizes the accuracy results, and the trajectories of body angular velocities (i.e., $\Omega(t)$) are illustrated in Figure~\ref{fig:rigid_body_traj}. Both feedback integrators under the adaptive gain selection and fixed gain $\alpha=\frac{1}{hL}$ achieve accuracies several orders of magnitude better than the unity gain variant. The unity gain method diverges for $h>10^{-3}$, whereas $\alpha=\tfrac{1}{hL}$ and the adaptive scheme maintains bounded error across the entire tested range of $h$. The adaptive scheme is marginally less accurate than fixed gain $\alpha=\tfrac{1}{hL}$, primarily because $\left\|\nabla^{2}V(x_k)\right\|$ varies little along the trajectory ($\left\|\nabla^{2}V(x_k)\right\|_F$ spans $[2334.63, 2412.56]$); with a safety factor $c=1.1$ and overapproximation with Frobenius norm, the local bound $c\left\|\nabla^{2}V(x_k)\right\|_F$ slightly overestimates the true maximum of $\left\|\nabla^{2}V\right\|$, yielding a more conservative gain and a small loss in performance. Notably, the adaptive gain feedback integrator requires no a priori estimate of the Lipschitz constant and achieves markedly higher accuracy than the unity gain feedback integrator, with only a marginal increase in computational cost.

Over the tested range of $h$, the splitting method incurs the highest computational cost but attains higher accuracy than feedback--integrator variants for $h > 10^{-6}$. In exact arithmetic, each partial flow is Hamiltonian and their symmetric composition is symplectic and time-reversible, which suppresses secular drift by making the scheme conserve a modified Hamiltonian~\cite{Hairer2006}. By contrast, feedback integrator applies a corrective term only after a deviation is present, so it cannot eliminate local truncation error related to local Hamiltonian flow a priori. 

For small step sizes, particularly for $h<10^{-4}$ in this example, the splitting method exhibits an upturn in the measured error. Over a fixed integration interval, decreasing $h$ increases the number of steps and hence the number of floating-point operations. Once the truncation error has been sufficiently reduced, it is well-known that accumulated round-off error can dominate the overall error. The observed upturn is therefore consistent with the standard round-off dominated regime. In contrast, the feedback-integrator variants do not show a comparable degradation over the tested range. A possible explanation is that the feedback term continuously corrects any deviations from the target set: since such correction is applied regardless of the cause of the deviation, whether it arises from truncation or round-off, it may also reduce the effect of accumulated floating-point perturbations in these experiments. Specialized finite-precision variants of symplectic schemes~\cite{rein2015whfast,earn2006symplectic} can mitigate such effects, but these approaches are problem- or integrator-specific and require nontrivial modifications. A theoretical analysis of round-off effects and finite-precision comparison are outside the scope of this work.
\subsection{The Kepler Problem}
We consider the following two--body dynamics in the Kepler problem.
\begin{subequations}\label{eq:dynamics_kepler}
\begin{equation} 
	\dot{x} = v,
\end{equation}
\begin{equation} 
	\dot{v} = -\mu \frac{x}{|x|^3},
\end{equation}
\end{subequations}
where $x \in \mathbb{R}^3 \setminus \{ (0,0,0) \}$ represents the position and $v \in \mathbb{R}^3$ the velocity. The gravitational parameter is denoted by $\mu$, where we use $\mu = 1$ throughout the simulations. Here, the angular momentum and Laplace--Runge--Lenz vector are the first integrals of the system, defined as follows respectively:
\begin{subequations}\label{eq:kepler_first_integrals}
\begin{equation} 
	L(x, v) = x \times v,
\end{equation}
\begin{equation} 
	A(x, v) = v \times (x \times v) - \mu\frac{x}{|x|}.
\end{equation}
\end{subequations}
The Lyapunov function is defined as 
\begin{equation}
	V(x, v) = \frac{k_1}{2}\left| L(x, v) - L_0 \right|^2 + \frac{k_2}{2}\left| A(x, v) - A_0 \right|^2,
\end{equation}
for $k_1 = 4$, $k_2 = 2$ where $L_0$ and $A_0$ represent the respective first integrals at initial value. We use $x_I = (1, 0, 0)$ and $v_I = (0, \sqrt{1.8}, 0)$ for the initial values, which correspond to $L_0 = (0, 0, \sqrt{1.8})$ and $A_0 = (0.8, 0, 0)$. The eccentricity of the Kepler orbit with initial value $(x_I, v_I)$ is $e = 0.8$, which corresponds to a period of $T = 70.2481$. We carry out integration on $[0, 1000T]$ interval. We estimate the Lipschitz constant of $\nabla V$ as $L \approx 515.4$. The gain update period for adaptive feedback integrator is set as $T_{\mathrm{update}} = 0.1$ seconds. We obtain the results for $h \in \left[10^{-6}, 10^{-1}\right]$, with comparison with a benchmark framework, the \emph{Störmer--Verlet method}~\cite{Hairer2006}.

\begin{figure}[t]
  \centering
  \begin{subfigure}[t]{1.0\linewidth}
    \includegraphics[width=\linewidth]{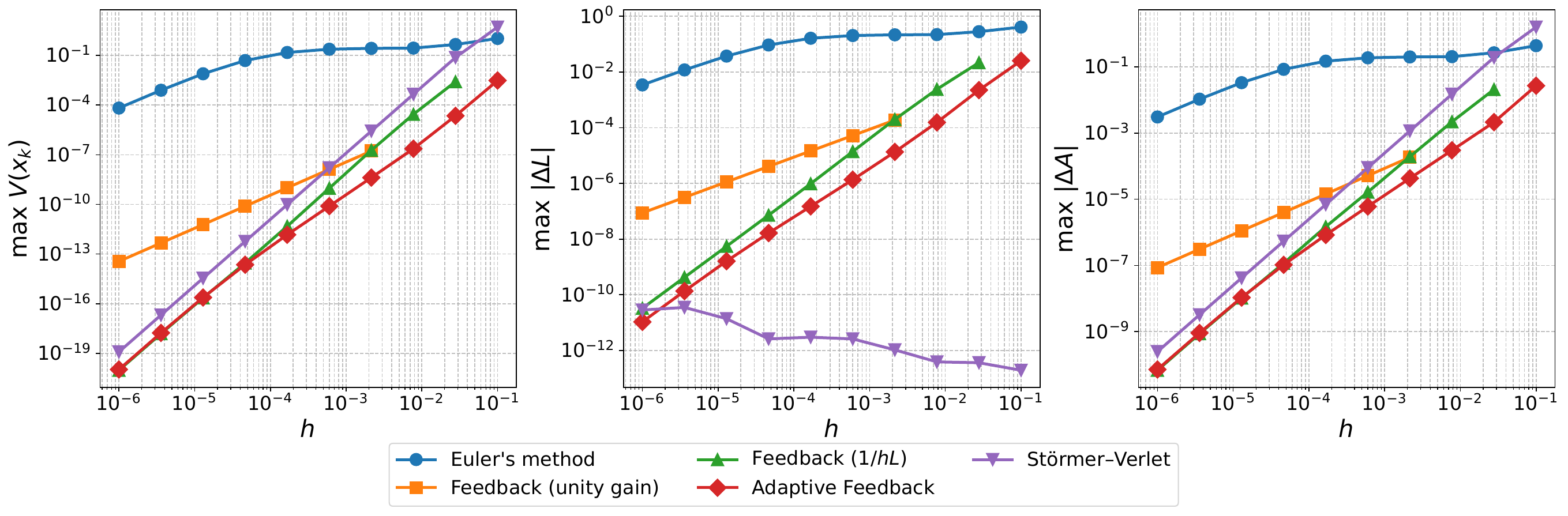}
    \caption{Accuracy vs.\ $h$}
  \end{subfigure}\hfill
  \begin{subfigure}[t]{0.5\linewidth}
    \includegraphics[width=\linewidth]{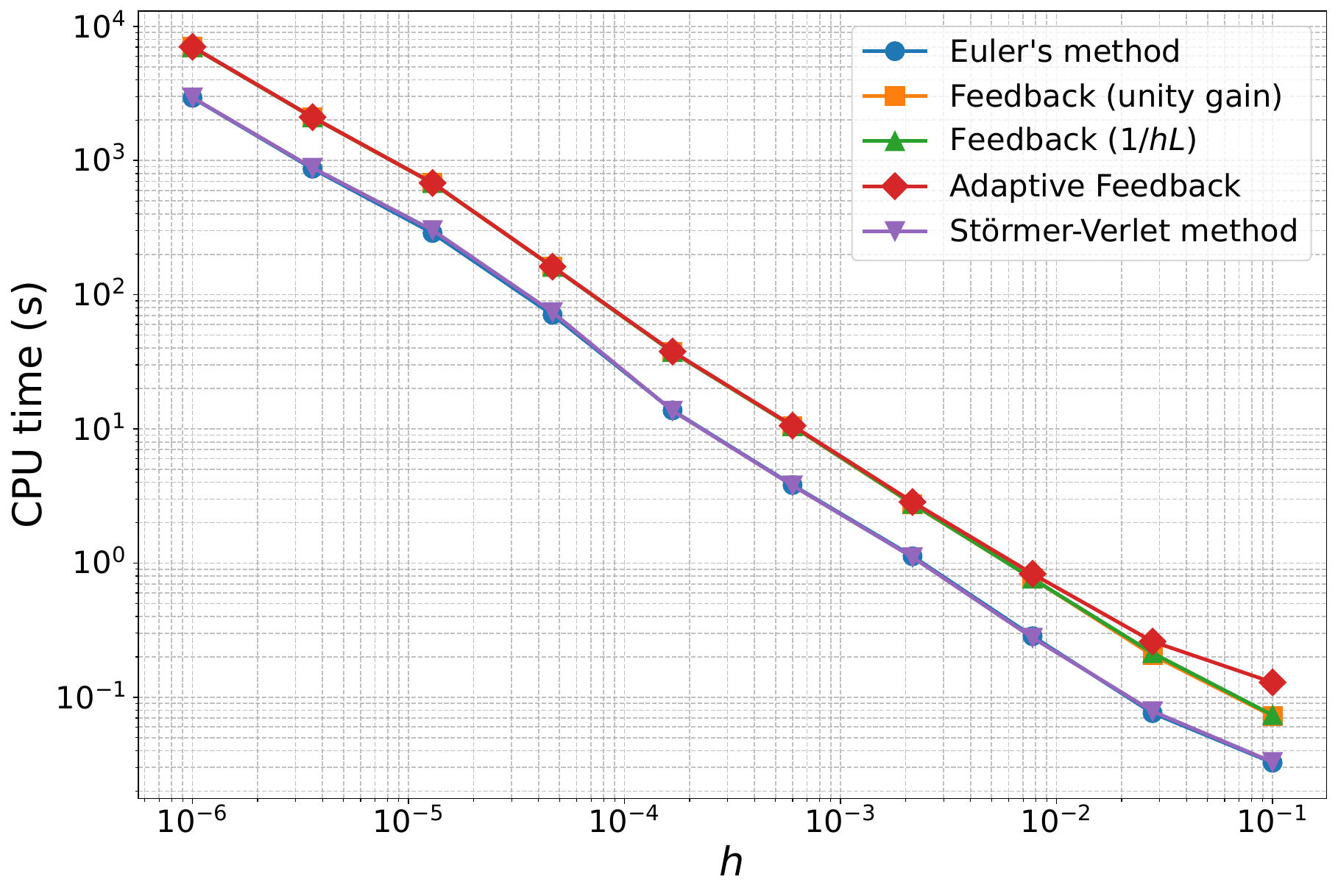}
    \caption{CPU time vs.\ $h$}
  \end{subfigure}
  \caption[Accuracy and cost vs.\ step size $h$]%
  {\textbf{Accuracy results of the Kepler problem.} Integration with Euler's method over $1000$ periods with $T = 70.2481$. (a) maximum $V(x_k)$ and deviation of first integrals along the trajectories. (b) CPU time dedicated for each integration scheme.} \label{fig:kepler_results}
\end{figure}

\begin{figure}[t]
  \centering
  \begin{subfigure}[t]{0.32\linewidth}
    \includegraphics[width=\linewidth]{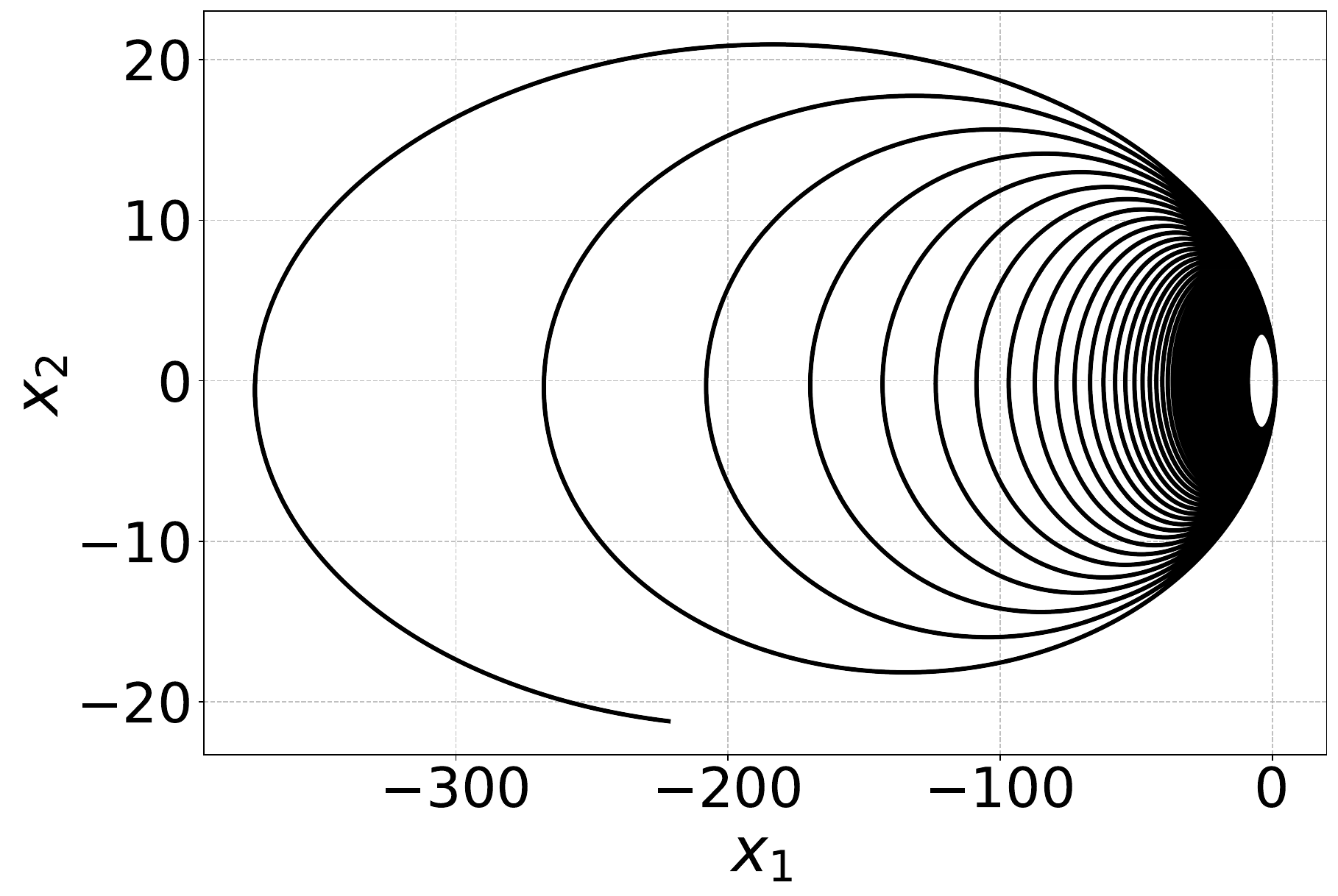}
    \caption{Euler's method}
  \end{subfigure}\hfill
  \begin{subfigure}[t]{0.32\linewidth}
    \includegraphics[width=\linewidth]{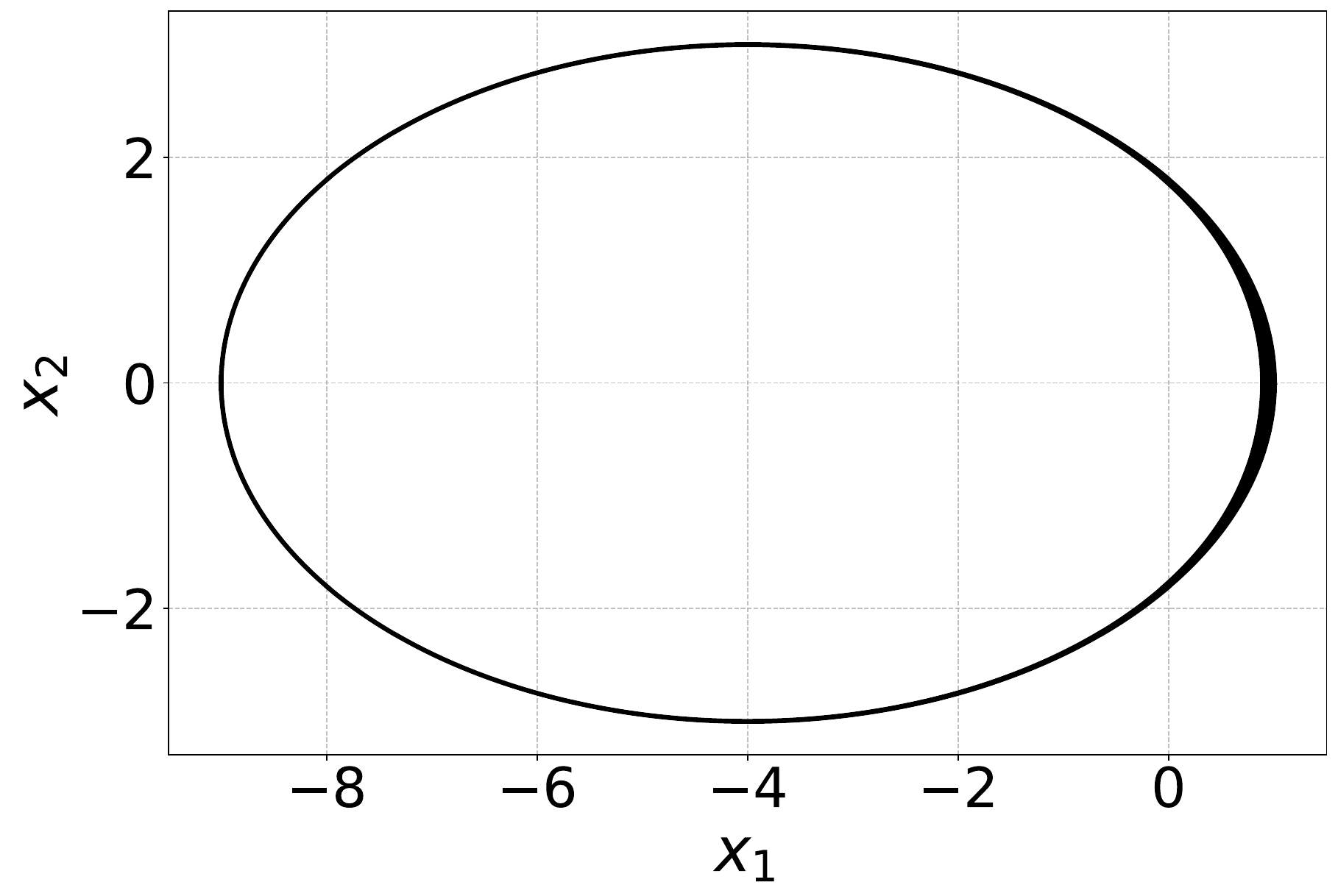}
    \caption{Feedback ($\alpha = 1$)}
  \end{subfigure}\hfill
  \begin{subfigure}[t]{0.32\linewidth}
    \includegraphics[width=\linewidth]{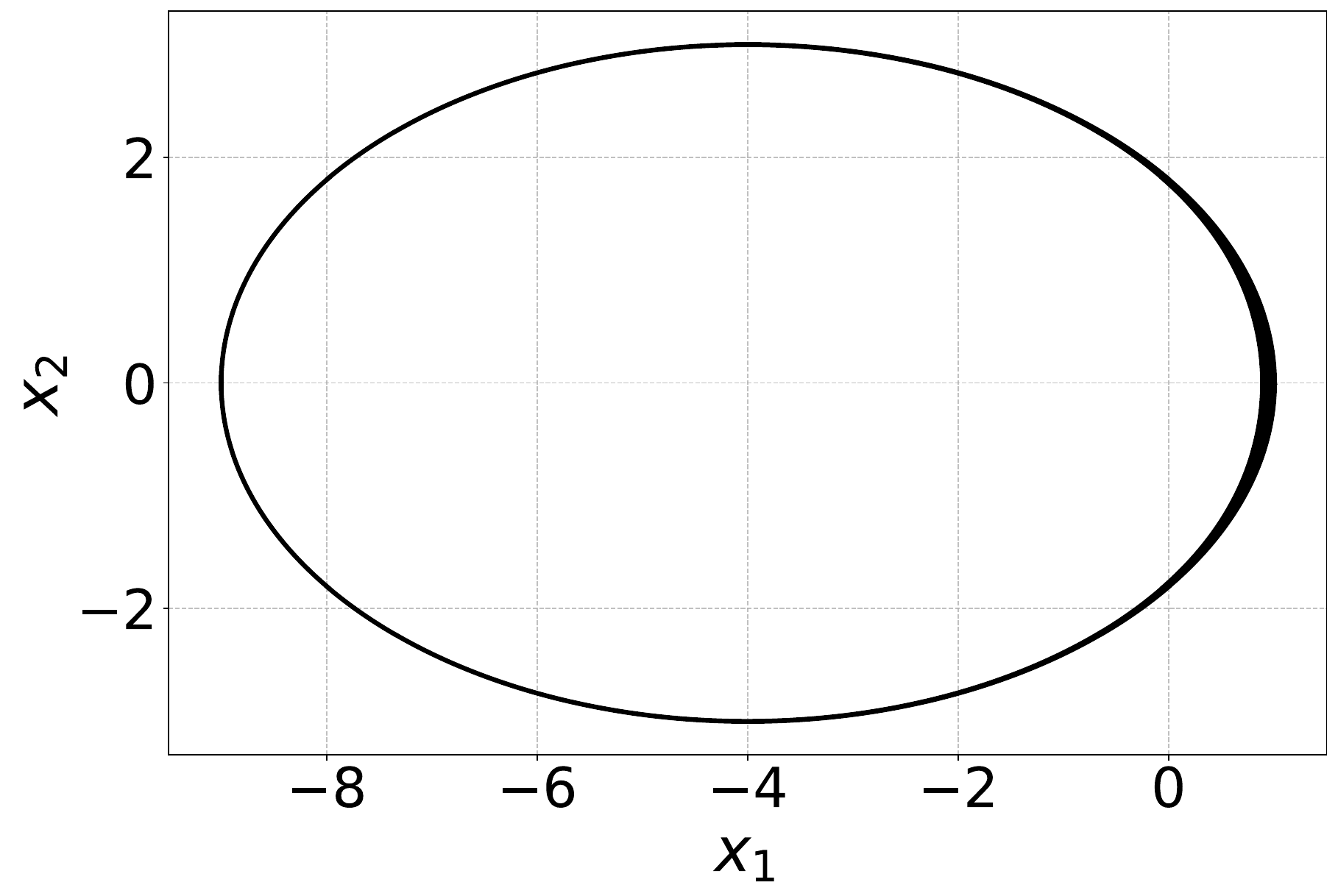}
    \caption{Feedback ($\alpha = 1/hL$)}
  \end{subfigure}\hfill
  \begin{subfigure}[t]{0.48\linewidth}
    \includegraphics[width=\linewidth]{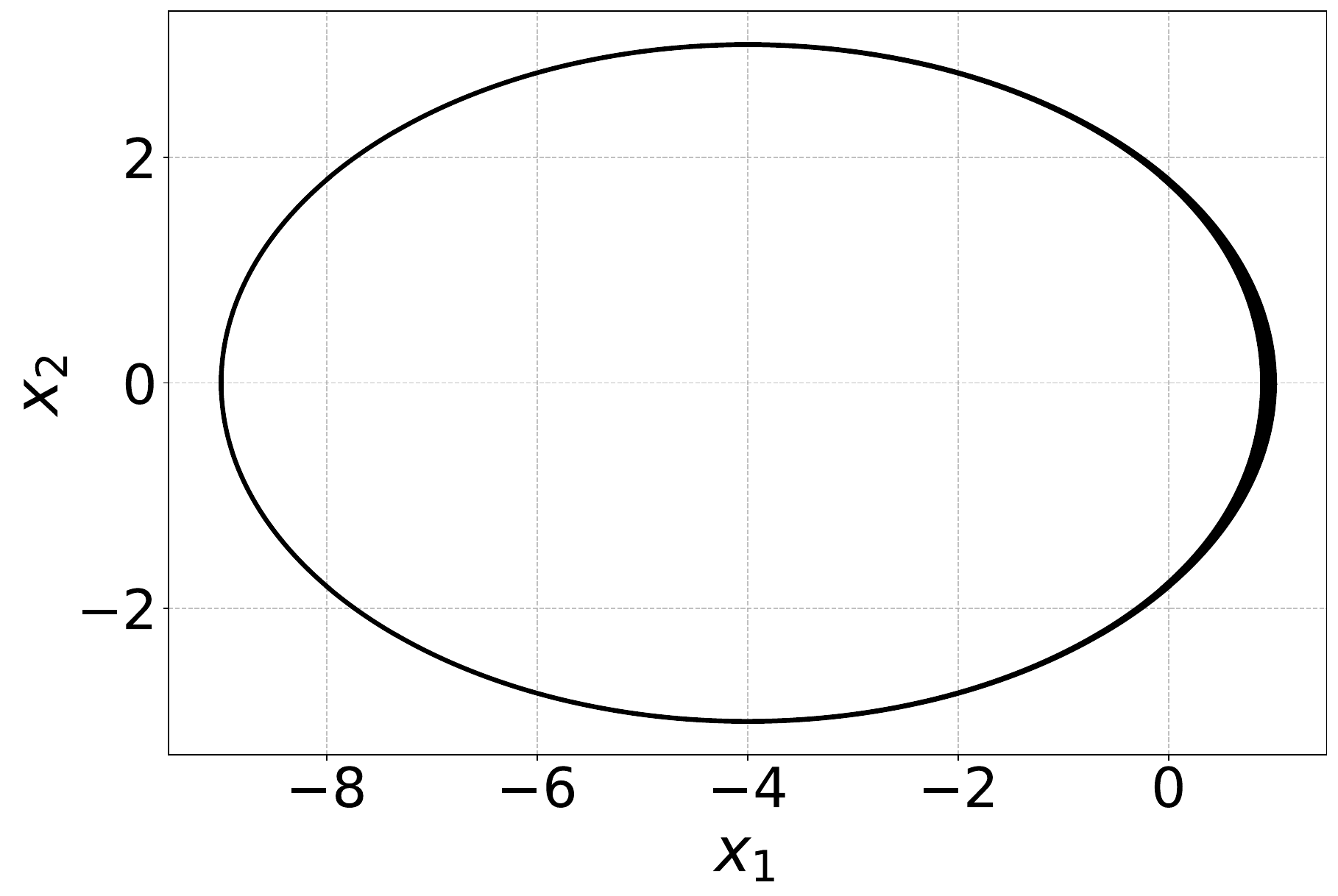}
    \caption{Adaptive feedback}
  \end{subfigure}\hfill
  \begin{subfigure}[t]{0.48\linewidth}
    \includegraphics[width=\linewidth]{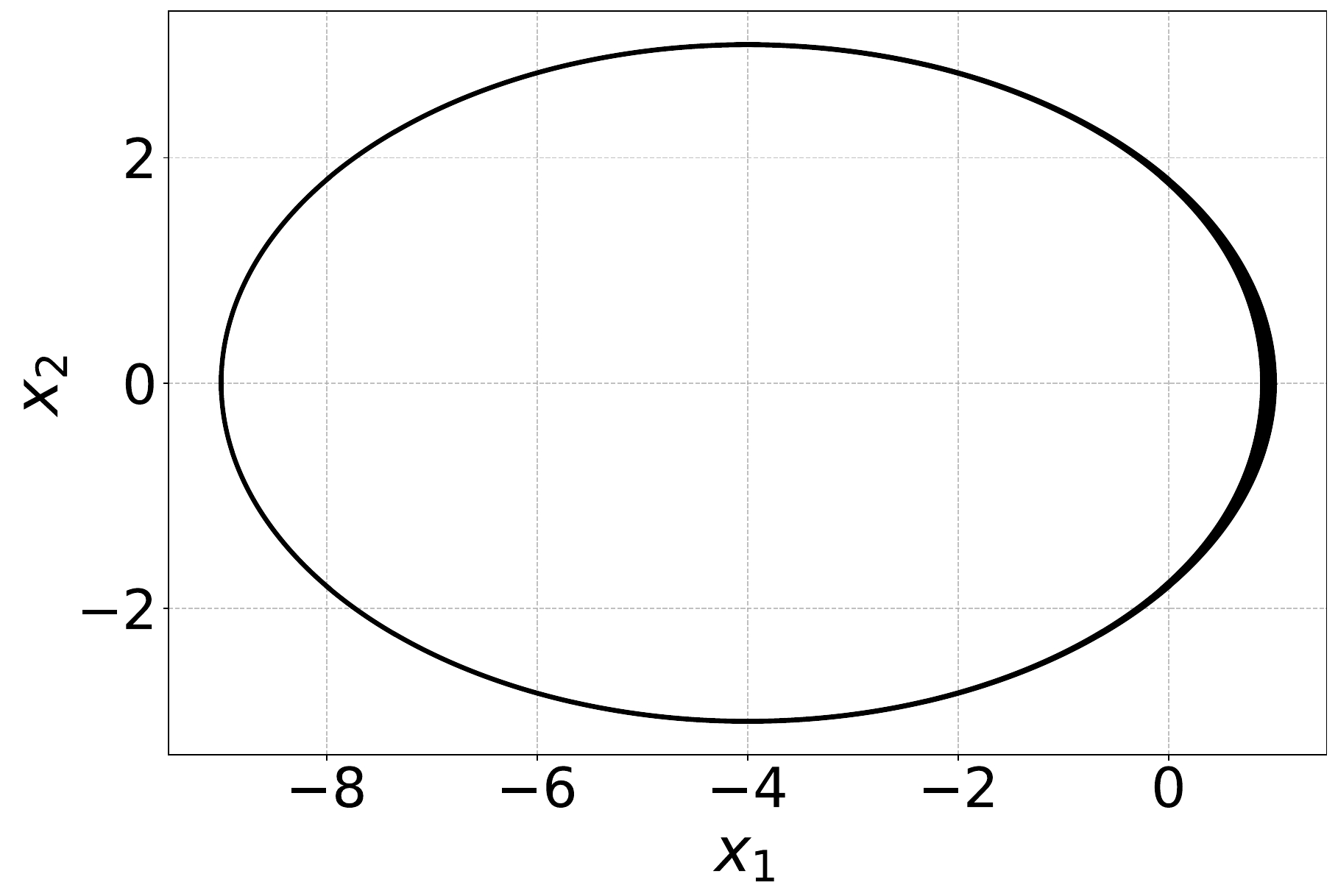}
    \caption{Störmer-Verlet}
  \end{subfigure}\hfill
  \caption[Kepler Trajectories]
  {\textbf{Trajectories of the Kepler Problem.} Integration over $1000$ periods with $T = 70.2481$ and $h = 10^{-3}$. Feedback integrators are implemented with Euler's method.} \label{fig:kepler_traj}
\end{figure}

Figure~\ref{fig:kepler_results} presents the accuracy results and Figure~\ref{fig:kepler_traj} illustrates the trajectories. Over the entire tested range of $h$, the adaptive gain feedback integrator keeps the error bounded. The Taylor-based fixed gain rule $\alpha=\tfrac{1}{hL}$ diverges only at the largest step $h = 10^{-1}$, whereas the unity gain scheme diverges for $h>10^{-8/3}\approx 2.15\times10^{-3}$. Relative to unity gain, these two variants achieve errors lower by several orders of magnitude. In this problem the adaptive scheme is more accurate than the fixed gain $\alpha=\tfrac{1}{hL}$, which is consistent with the pronounced variability of $\left\|\nabla^2 V(x_k)\right\|$ that spans $[30.22, 515.4]$ (about $17\times$). Thus, periodic gain update better tracks the changing Lipschitz scale than a single global setting. Except at the largest step size where the fixed gain rule $\alpha=\tfrac{1}{hL}$ diverges, feedback integrators with fixed and adaptive gain selection achieve the lowest errors among the tested schemes, at the expense of only a slight increase in computational cost. For these feedback integrator variants, runtime is dominated by evaluations of $\nabla V$, so the cost depends strongly on the chosen $V$. Both feedback schemes also outperform the Störmer-Verlet method in accuracy at comparable step sizes, albeit with modestly higher cost.
\subsection{Perturbed Kepler Problem}
For the last example, we consider a perturbed Kepler problem with rotational symmetry. The dynamics is as follows,
\begin{subequations}\label{eq:dynamics_perturbed_kepler}
\begin{equation} 
	\dot{x} = v,
\end{equation}
\begin{equation} 
	\dot{v} = -U'\left(|x|\right) \frac{x}{|x|},
\end{equation}
\end{subequations}
with $x \in \mathbb{R}^3 \setminus \{ (0,0,0) \}$ representing the position and $v \in \mathbb{R}^3$ the velocity. The potential function $U: (0, \infty) \rightarrow \mathbb{R}$ is assumed to be radially symmetric. The first integrals of this dynamics are total mechanical energy and angular momentum, defined as follows:
\begin{subequations}\label{eq:perturbed_kepler_first_integrals}
\begin{equation} 
	E(x, v) = \frac{1}{2}|v|^2 + U\left(|x|\right),
\end{equation}
\begin{equation} 
	L(x, v) = x \times v.
\end{equation}
\end{subequations}
Accordingly, the Lyapunov function is defined as 
\begin{equation}
	V(x, v) = \frac{k_1}{2}\left| E(x, v) - E_0 \right|^2 + \frac{k_2}{2}\left| L(x, v) - L_0 \right|^2
\end{equation}
for constants $k_1 = 2$, $k_2 = 3$, and $E_0$ and $L_0$ represent the respective first integrals at initial value. Throughout the demonstration, we use $U\left(|x|\right) = -\frac{\mu}{|x|} - \frac{\delta}{|x|^3}$ with $\mu = 1$ and $\delta = 0.0025$. For eccentricity $e = 0.6$, we take $x_I = (1-e, 0, 0) = (0.4, 0, 0)$ and $v_I = \left(0,\sqrt{\frac{1+e}{1-e}},0\right) = (0, 2, 0)$, yielding $E_0 \approx -0.5391$ and $L_0 = (0, 0, 0.8)$. We integrate over $[0, 200]$. We estimate the Lipschitz constant of $\nabla V$ as $L \approx 148.03$, and along the trajectory $\|\nabla^2 V\|$ lies in the range $[6.27, 148.03]$. For adaptive gain scheme, we set the gain update period as $T_{\mathrm{update}} = 0.1$. We report results for $h \in \left[ 10^{-9}, 10^{-1} \right]$, and include the Störmer-Verlet method as a benchmark.

\begin{figure}[t]
  \centering
  \begin{subfigure}[t]{1.0\linewidth}
    \includegraphics[width=\linewidth]{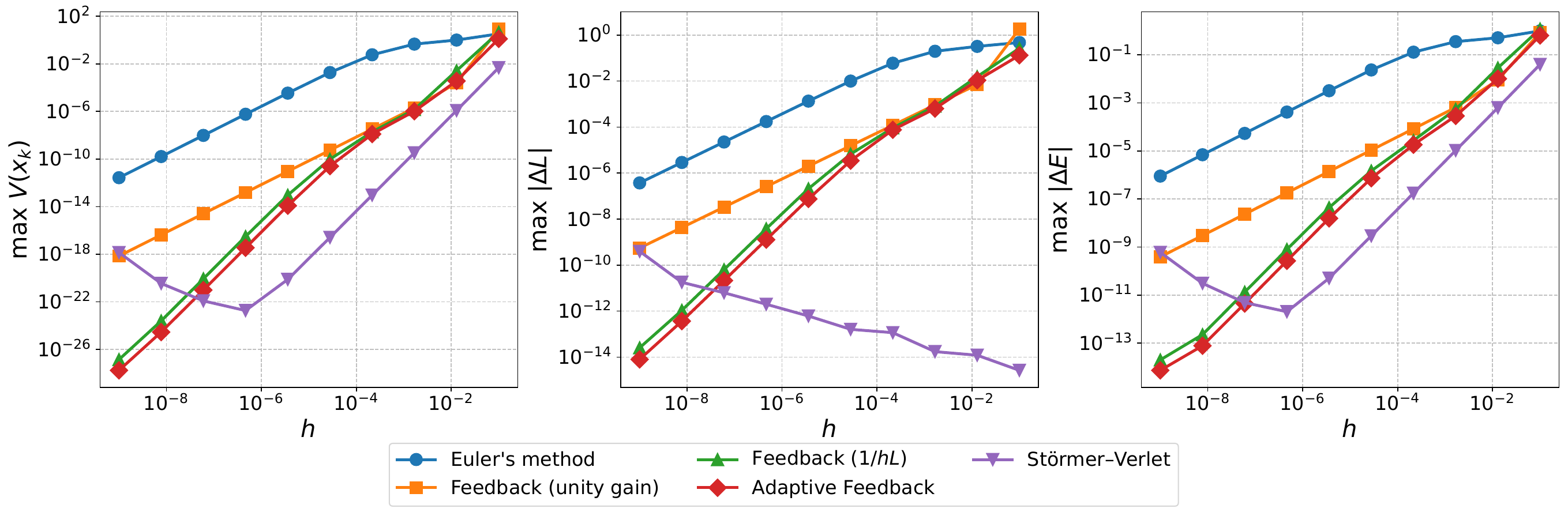}
    \caption{Accuracy vs.\ $h$}
  \end{subfigure}\hfill
  \begin{subfigure}[t]{0.5\linewidth}
    \includegraphics[width=\linewidth]{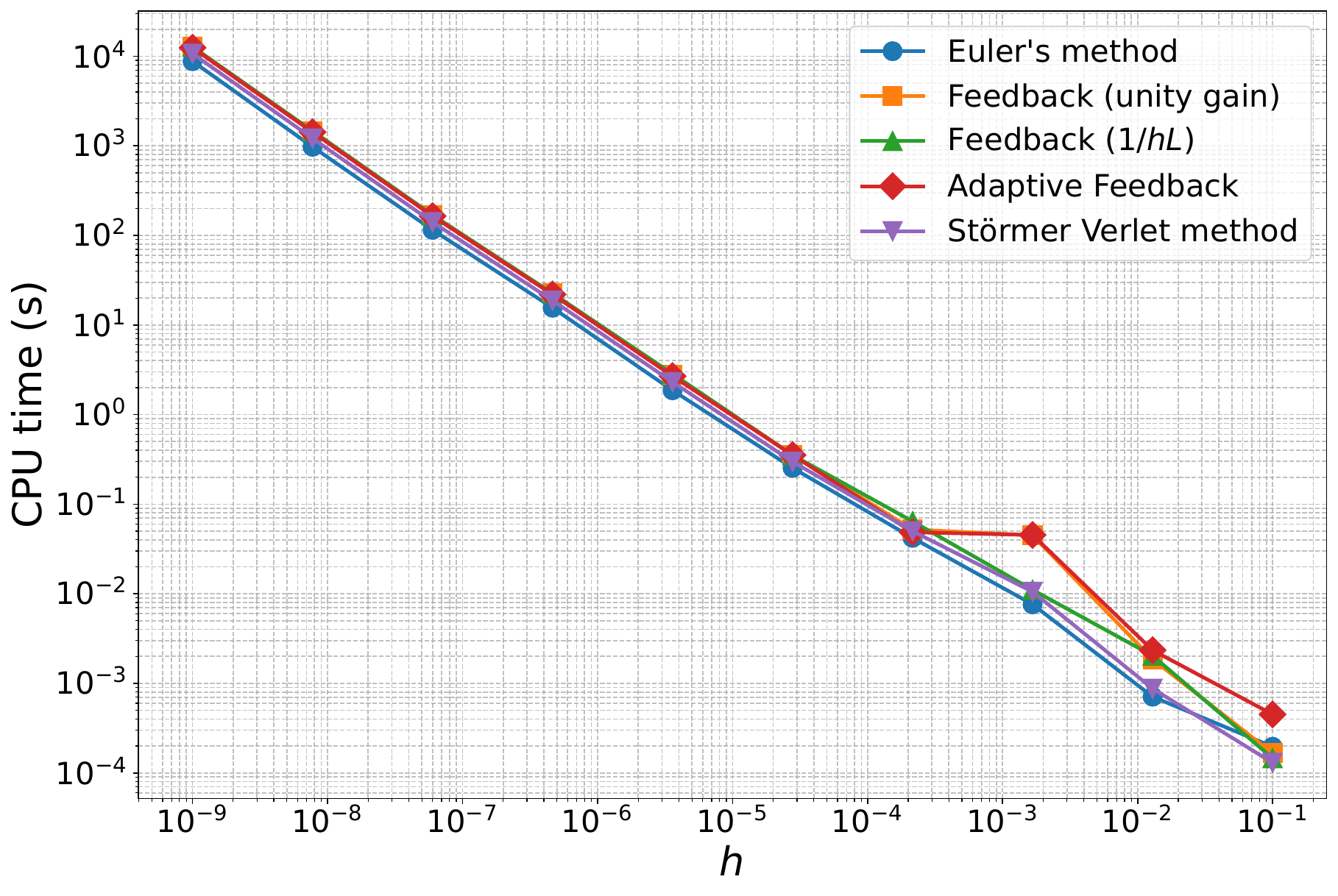}
    \caption{CPU time vs.\ $h$}
  \end{subfigure}
  \caption[Accuracy and cost vs.\ step size $h$]%
  {\textbf{Accuracy results of perturbed Kepler problem.} Integration with Euler's method over $[0, 200]$. (a) maximum $V(x_k)$ and deviation of first integrals along the trajectories. (b) CPU time dedicated for each integration scheme.} \label{fig:perturbed_kepler_results}
\end{figure}

\begin{figure}[t]
  \centering
  \begin{subfigure}[t]{0.32\linewidth}
    \includegraphics[width=\linewidth]{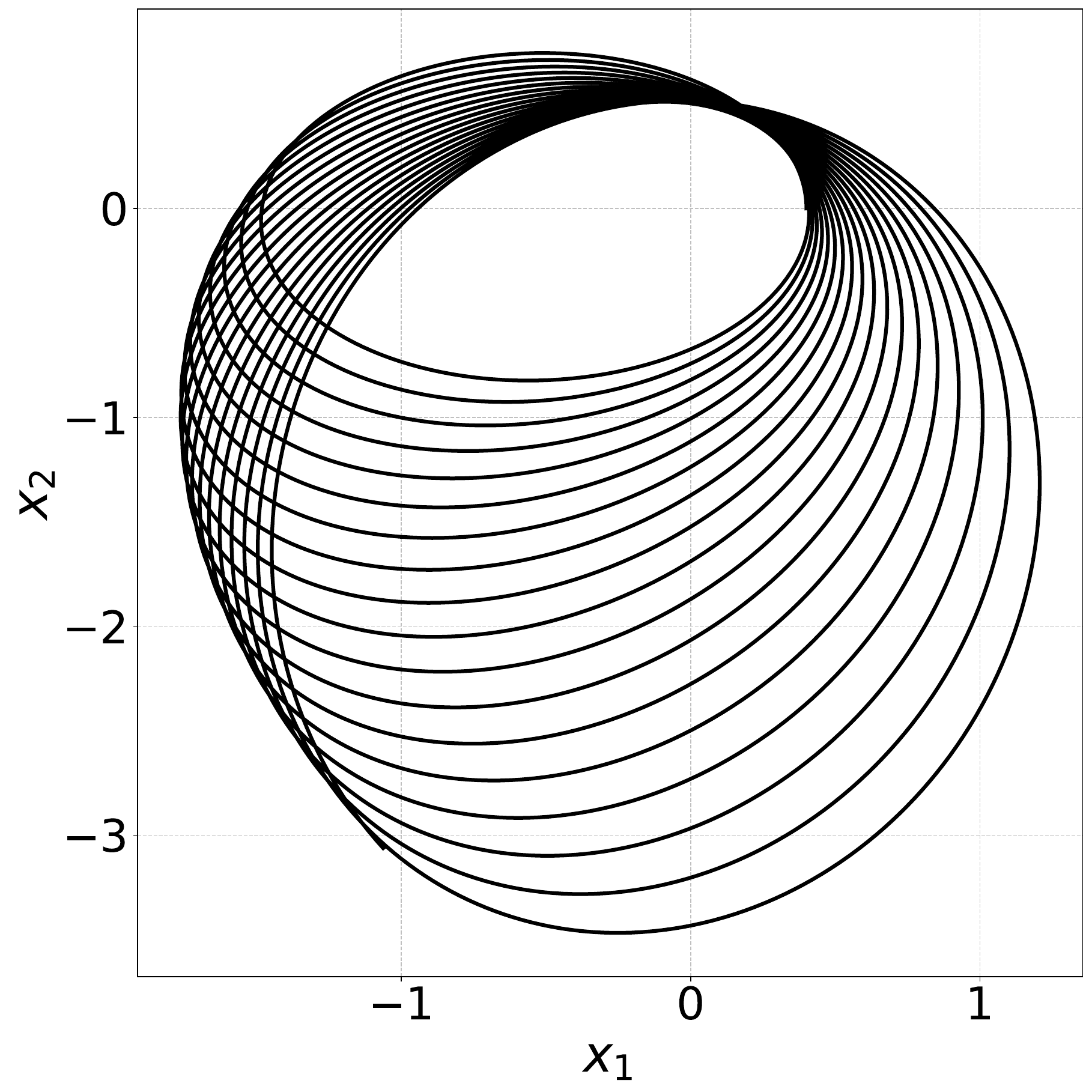}
    \caption{Euler's method}
  \end{subfigure}\hfill
  \begin{subfigure}[t]{0.32\linewidth}
    \includegraphics[width=\linewidth]{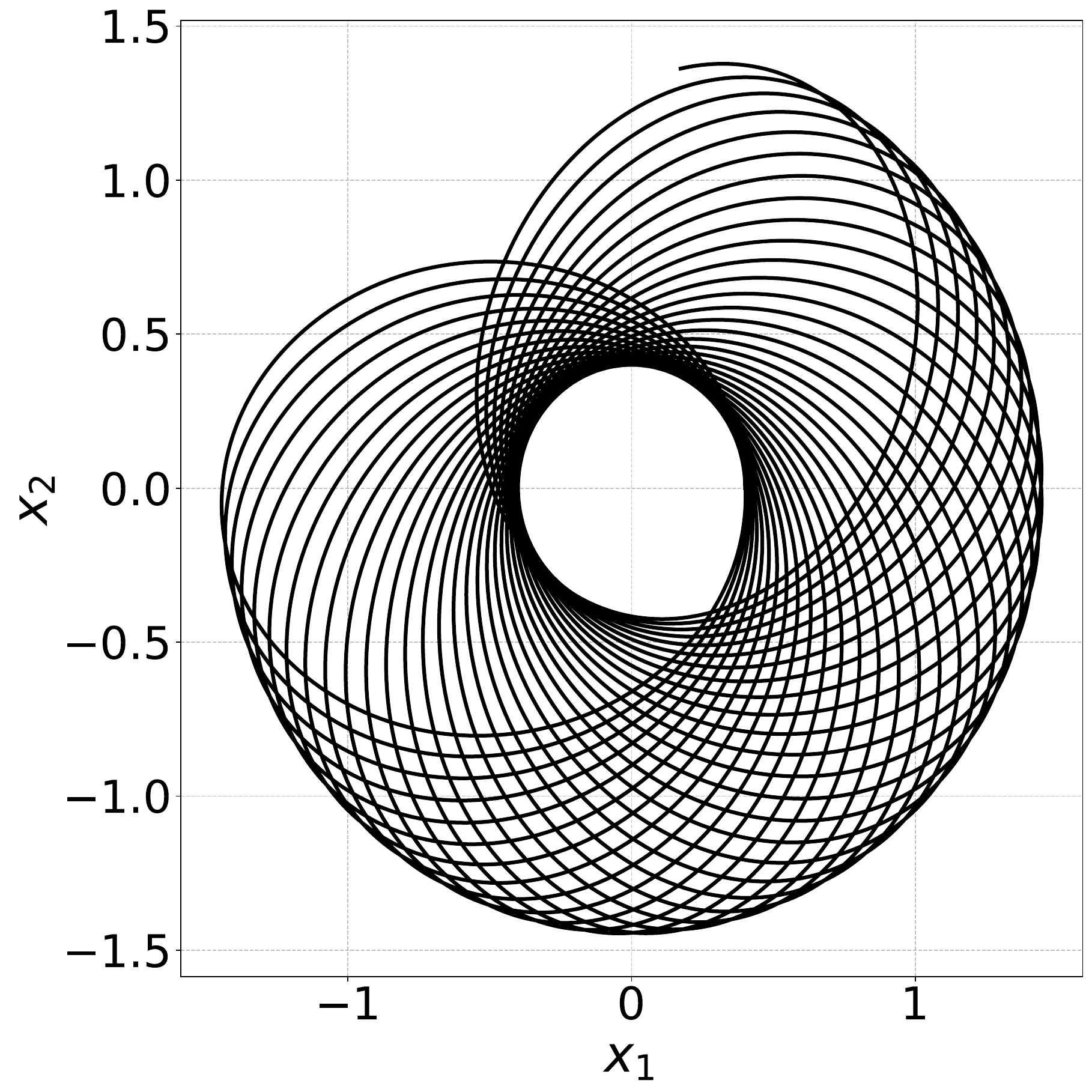}
    \caption{Feedback ($\alpha = 1$)}
  \end{subfigure}\hfill
  \begin{subfigure}[t]{0.32\linewidth}
    \includegraphics[width=\linewidth]{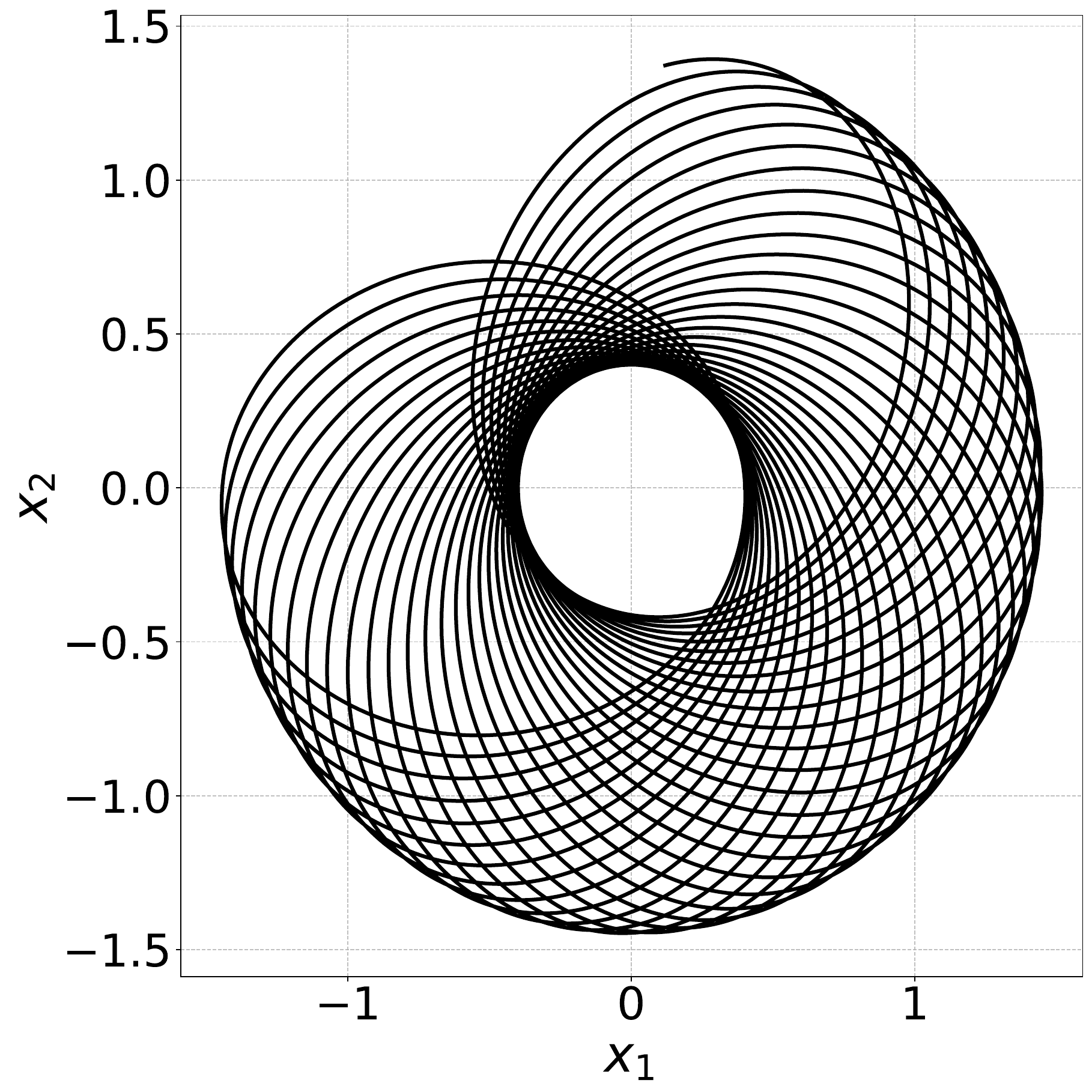}
    \caption{Feedback ($\alpha = 1/hL$)}
  \end{subfigure}\hfill
  \begin{subfigure}[t]{0.48\linewidth}
    \includegraphics[width=\linewidth]{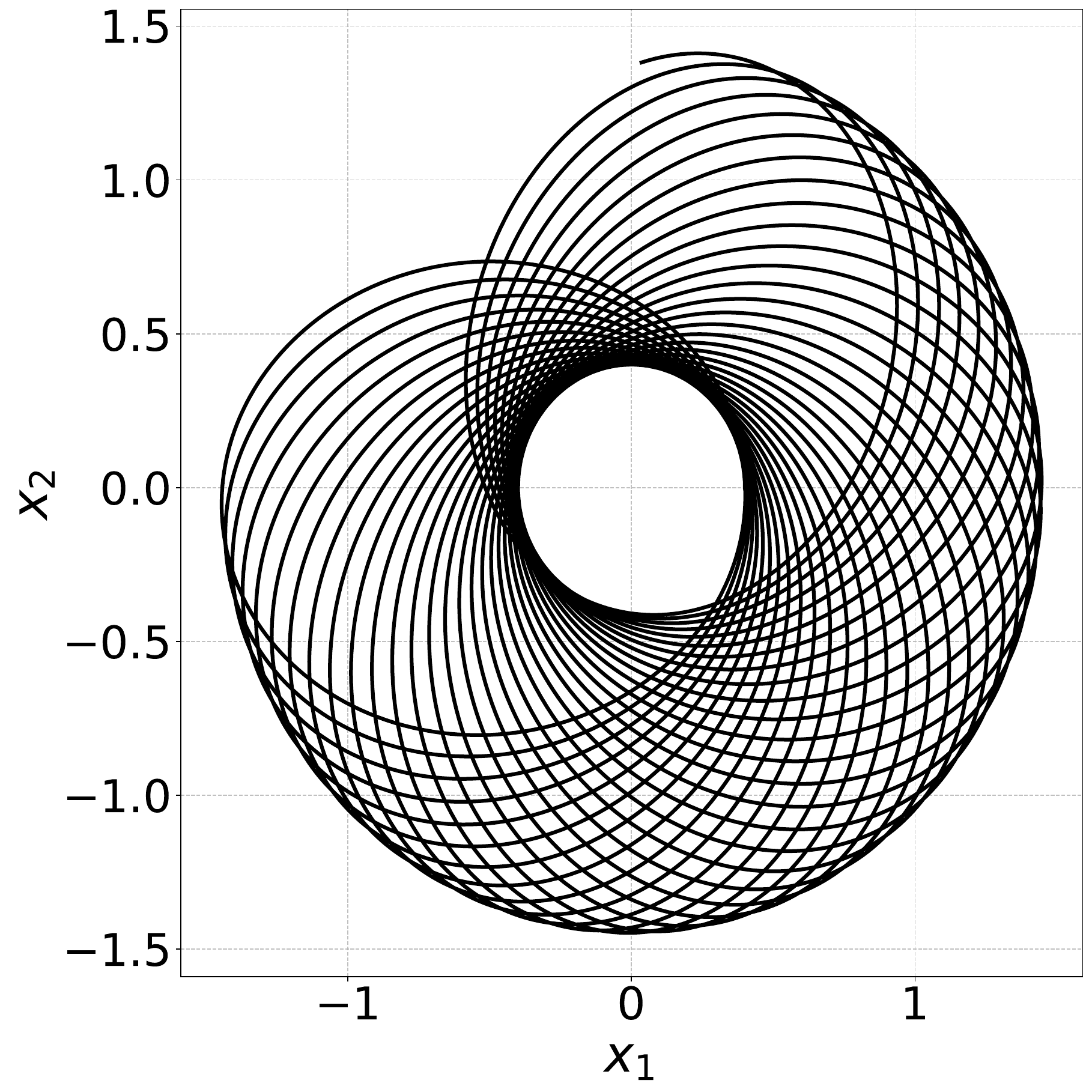}
    \caption{Adaptive feedback}
  \end{subfigure}\hfill
  \begin{subfigure}[t]{0.48\linewidth}
    \includegraphics[width=\linewidth]{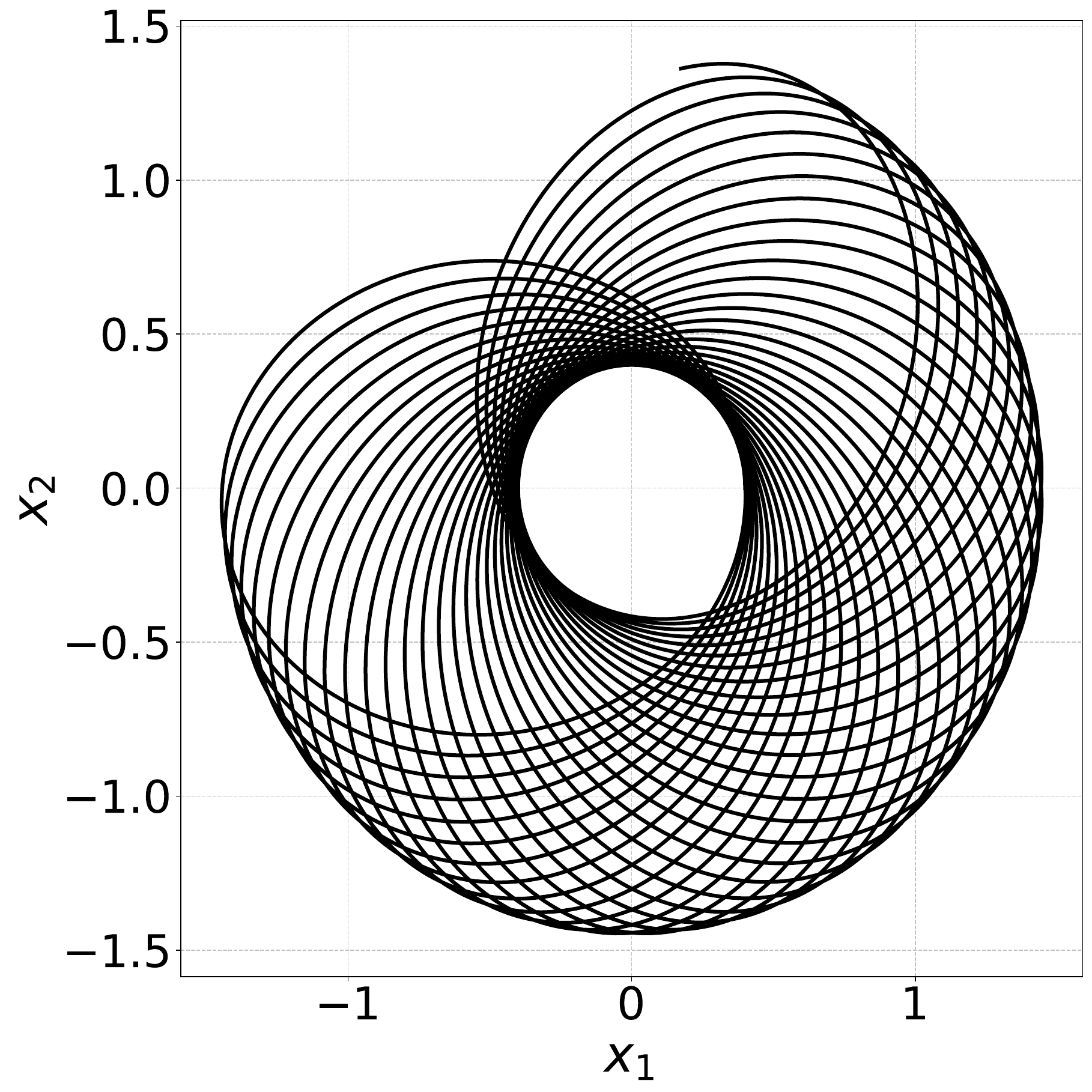}
    \caption{Störmer-Verlet}
  \end{subfigure}\hfill
  \caption[Perturbed Kepler Trajectories]
  {\textbf{Trajectories of perturbed Kepler Problem.} Integration on $[0, 200]$ with $h = 10^{-3}$. Feedback integrators are implemented with Euler's method.} \label{fig:perturbed_kepler_traj}
\end{figure}

Figure~\ref{fig:perturbed_kepler_results} presents the accuracy results and Figure~\ref{fig:perturbed_kepler_traj} illustrates the trajectories. Over the entire tested range of $h$, all feedback-integrator variants keep the error bounded. As in the Kepler example, the adaptive and the fixed gain rule $\alpha = \frac{1}{hL}$ outperform the unity gain variant by several orders of magnitude, with the adaptive scheme slightly more accurate. This is consistent with the variability of $\|\nabla^2 V\|$ spanning $[6.27, 148.03]$. The computational costs of all methods, including Störmer-Verlet, are comparable, and Störmer-Verlet yields the best accuracy at moderate step sizes $h \geq 10^{-65/9} \approx 5.995 \times 10^{-8}$. As $h$ is decreased below $10^{-6}$, the measured error of Störmer-Verlet exhibits an upturn, consistent with a round-off dominated regime, and eventually exceeds that of the feedback-integrator variants.
%
\section{Conclusion} \label{sec:conclusion}
This paper establishes a non-asymptotic foundation for feedback integrators and develops a complete small-step gain-selection theory in the Euler discretization setting. First, for general one-step discretizations, we proved that arbitrarily small sublevel neighborhoods of the target set are positively invariant for sufficiently small step sizes. This fills the gap in the original theory by providing a non-asymptotic-in-time bound from the first step onward. Second, under Euler discretization, we analyzed the scaled gain $\beta=h\alpha$. We characterized a range of scaled gains that guarantee positive invariance in the small-step regime and identified the scaling that minimizes the Taylor-based upper bound. We also developed stepwise and periodically updated adaptive gain-selection rules and proved corresponding boundedness guarantees for the resulting discrete trajectories.

Numerical demonstrations on free rigid body motion in $\operatorname{SO}(3)$, the Kepler problem, and a perturbed Kepler problem with rotational symmetry support the analysis. In the tested examples, the proposed gain-selection rules substantially reduce the error relative to the unity-gain baseline, with costs comparable to the baseline feedback integrator. The adaptive variants are particularly useful when the Hessian scale of the feedback Lyapunov function varies significantly along the trajectory, while the fixed rule can be competitive when this scale is nearly constant. In the smallest step-size regimes tested, some structure-preserving benchmarks exhibit an upturn in the measured error, consistent with accumulated round-off effects, whereas the feedback-integrator variants do not show comparable degradation in these experiments.

The results also clarify the role of Euler discretization in the feedback integrator framework. For general one-step methods, the non-asymptotic invariance result holds at the level of positive invariance, but the gain dependence is entangled with method-dependent local error terms. In contrast, under Euler discretization the method-dependent remainder terms vanish, and the scaled gain enters the one-step bound explicitly. Thus, Euler discretization provides the canonical explicit setting in which gain selection admits a closed-form small-step theory: scaled gains guaranteeing positive invariance can be characterized, the Taylor-bound-minimizing scaling can be identified, and adaptive variants can be equipped with boundedness guarantees.

Extending the gain selection theory beyond the Euler discretization setting remains an important open direction. Higher-order one-step methods introduce method-dependent mixed terms in step sizes and feedback gain, so comparable gain-selection rules would have to account for the specific discretization. A separate direction is to develop a finite-precision analysis of feedback integrators and to understand whether the empirical round-off behavior observed here can be explained theoretically.
%
%
%
\section*{Acknowledgement} \label{sec:ack}
Most parts of this work were done when the first author was an undergraduate student at the School of Electrical Engineering, KAIST. This work was supported by the National Research Foundation of Korea(NRF) grant funded by the Korea government(MSIT) (RS-2026-25473622).
%


%
\bibliography{Ref_Papers.bib}

\end{document}